\documentclass[reqno,12pt]{birkjour}
\usepackage{graphicx, array}
\usepackage{amssymb, amsmath, amsthm, enumerate, color, 
soul, longtable, nicefrac, xtab}
\usepackage{tablefootnote}
\usepackage[parfill]{parskip}    
\usepackage{tikz}
\usetikzlibrary{calc, positioning, through, math}
\usetikzlibrary{backgrounds}
\usepackage{tabularx}

\usepackage{caption}
\usepackage{subcaption}

\usepackage[]{floatrow}
\usepackage[color links=true, backref=page]{hyperref}

\tikzstyle{vtx}=[inner sep=1pt,draw, shape=circle]
\tikzstyle{line}=[inner sep=3pt,draw, shape=circle, fill=black!20]
\tikzset{>=stealth}

\newcommand{\darcl}[3]
{\draw[] (#1) edge [<-, bend left=20] node[midway, fill=white, inner sep=1 pt, font = \tiny]{$#2$} (#3);
\draw[ultra thick] (#3)edge [bend left=20]  (#1);}

\newcommand{\darcr}[3]
{\draw (#1) edge [<-, bend right=20] node[midway, fill=white, inner sep=1 pt,font = \tiny]{$#2$} (#3); 
\draw[ultra thick] (#3)edge [bend right=20]  (#1);}

\definecolor{mycyan}{RGB}{0, 255, 255}
\definecolor{mymagenta}{RGB}{255,0, 255}

\newcommand{\du}{\ensuremath{\mathbf{DU}}}

\newcommand{\AR}{\ensuremath{\mathbf{AR}}}

\newcommand{\PS}{\ensuremath{\mathbf{PS}}}

\newcommand{\AS}{\ensuremath{\mathbf{AS}}}

\DeclareGraphicsExtensions{.pdf}

 \newtheorem{thm}{Theorem}[section]
 \newtheorem{cor}[thm]{Corollary}
 \newtheorem{lem}[thm]{Lemma}
 \newtheorem{prop}[thm]{Proposition}
 \theoremstyle{definition}
 \newtheorem{defn}[thm]{Definition}
 \theoremstyle{remark}
 \newtheorem{rem}[thm]{Remark}
 
 \numberwithin{equation}{section}
  \newtheorem*{construction}{Construction}

  \newtheorem{question}{Question}

\newcommand{\be}{\begin{enumerate}}
\newcommand{\ee}{\end{enumerate}}

\newcommand{\bi}{\begin{itemize}}
\newcommand{\ei}{\end{itemize}}

\newcommand{\mc}[1]{\mathcal{#1}}

\newenvironment{matr}[1]{\left( \begin{array}{#1}}{\end{array}\right)}

\title[Bounds on $(n_{5})$ and $(n_{6})$ configurations]{ 
New bounds on the existence of $(n_{5})$ and $(n_{6})$ configurations: the Gr\"{u}nbaum Calculus revisited\\}
\author{Leah Wrenn Berman$^{*}$\footnotetext{*Corresponding Author}}
\address{Department of Mathematics \& Statistics \\ University of Alaska Fairbanks\\ P.O. Box 756660\\  Fairbanks, Alaska, 99775-6660\\ USA}
\email{lwberman@alaska.edu}
\author{G\'{a}bor G\'{e}vay}
\address{Bolyai Institute\\University of Szeged\\Aradi v\'ertan\'uk tere 1\\ Szeged, 6720 \\ Hungary}
\email{gevay@math.u-szeged.hu}
\author{Toma\v{z} Pisanski}
\address{University of Primorska\\ Koper, Slovenia\\ and
Institute of Mathematics, Physics and Mechanics\\University of Ljubljana\\ Ljubljana, Slovenia}
\email{pisanski@upr.si}
\date{\today}							

\begin{document}
\begin{abstract}
The ``Gr\"unbaum Incidence Calculus'' is the common name of a collection of operations introduced by Branko Gr\"unbaum 
to produce new $(n_{4})$ configurations from various input configurations.  In a previous paper, we generalized two of these 
operations to produce operations on arbitrary $(n_k)$ configurations, and we showed that for each $k$, there exists an integer 
$N_{k}$ such that for all $n \geq N_{k}$, there exists at least one $(n_{k})$ configuration, with current records $N_{5}\leq 576$ 
and $N_{6}\leq 7350$. In this paper, we further extend the Gr\"unbaum calculus; using these operations, as well as a collection of previously known and novel ad hoc constructions, we refine the bounds for $k = 5$ and $k = 6$.
Namely, we show that $N_5 \leq 166$ and $N_{6}\leq 585$. 
\end{abstract}

\maketitle

\section{Introduction}

A \emph{geometric $(n_{k})$ configuration} is a collection of $n$ points and $n$ straight lines in the Euclidean plane such that each line passes through $k$ 
points and each point lies on $k$ lines. In a series of papers \cite{Gru2000, Gru2000b,Gru2002,Gru2006} and in his 2009 book on configurations \cite{Gru2009b}, 
Branko Gr\"unbaum described a sequence of operations to produce new $(n_{4})$ configurations from various input configurations. These operations were later 
called the ``Gr\"unbaum Incidence Calculus'' (see, e.g., \cite[p. 251]{PisSer2013}). In \cite{BerGevPis2021}, we generalized two of those constructions, which 
we called \emph{affine replication} and \emph{affine switch}, to produce operations on arbitrary $(n_{k})$ configurations, and we used those two constructions 
to show that for any $k\geq 2$ there exists an integer $N_k$ such that for any $n, n \geq N_k$ there exists at least one geometric $(n_k)$ configuration. We 
refer to $(n_k)$ as the \emph{symbol} and $k$ as the \emph{type} of a configuration. Similarly, in case of a non-balanced  configuration (i.e.\ where the number 
of points and number of lines are different), we use $(p_q, n_k)$ for its symbol and $[q,k]$ for its type. 

In this paper, we generalize further the affine replication and affine switch operation, as well as Gr\"unbaum's two \emph{deleted union} constructions; we also
describe an additional operation that we call \emph{parallel switch}. Along with these construction techniques, we describe or introduce various ad hoc construction 
methods that produce either sporadic $(n_{5})$ (resp. $(n_{6})$) configurations or systematic constructions that produce arithmetic series of $(n_{5})$ (resp.\
$(n_{6})$) configurations. By combining ad hoc constructions and the old and new Gr\"unbaum Calculus operations, we decrease the bounds on $N_{5}$ and 
$N_{6}$ significantly: we decrease $N_{5}$ from $576$ to $166$, and $N_{6}$ from 7350 to $585$.

We note that Garrett Flowers, in his 2015 PhD thesis \cite{Flo2015}, used ideas from the Gr\"unbaum Calculus to show that $N_{k}$ always exists --- in fact, he 
shows a stronger result, that for any $r$, $k$, there exists an integer $N'(r,k)$ such that an $[r, k]$-configuration on $n$ points exists for all $n \geq N'(r, k)$ 
satisfying the divisibility condition $nr = bk$. However, his bounds are quite large: he says they are ``roughly on the order of $(k^{r})^{2}$'', so in the case 
we are interested in, where $r = k$, the bounds provided in \cite{BerGevPis2021} and refined in this paper are much smaller.

The structure of the paper is as follows. In Section \ref{sect:tools}, we describe and refine some geometric tools which allow us to carefully and specifically define 
the geometric constructions we will discuss in the next sections. In Section \ref{sect:IncCalc}, we describe the Gr\"unbaum Calculus operations we will be using. 
Section \ref{sect:n5systematic} describes previously discovered systematic constructions for 5-configurations and introduces several new constructions. In Section 
\ref{sect:bound5cfgs}, using the systematic and Gr\"unbaum calculus constructions, we prove that $N_{5} \leq 166$. In Section \ref{sect:bound6cfgs}, we review 
known constructions for 6-configurations and use those, in addition to the previous results, to show that $N_{6} \leq 585$.


\section{Some conceptual tools} \label{sect:tools}


\subsection{Flexible and compatible configurations} \label{subsect:flex}

Here we introduce two notions which are useful when applying the various operations on configurations discussed later on in this paper.

\begin{defn} \label{def:flex}
Let $\mathcal C$ be a $k$-configuration and $\ell$ one of its lines. We say that the line $\ell$ is \emph{flexible} if one can prescribe the 
position of $k$ configuration points on some line $\ell'$ of a configuration $\mathcal C'$, and there exists a projective transformation that 
maps $\mathcal C$ to  $\mathcal C'$ in such a way that $\ell$ gets mapped to $\ell'$ and the $k$ configuration points on $\ell$ are mapped 
into the prescribed positions on $\ell'$. Moreover, any configuration admitting a flexible line is called \emph{line-flexible}. Using duality one 
may define analogously a \emph{flexible point} and a \emph{point-flexible} configuration.
\end{defn}

\begin{defn} \label{def:comp}
Let $\mathcal C_1$ and $\mathcal C_2$ be two $k$-configurations. We say that $\mathcal C_1$ and $\mathcal C_2$ are \emph{compatible}
if  there is a point $P \in \mathcal C_1$ and a pencil $\mathcal P$ of lines $\ell_1\dots, \ell_k \in \mathcal C_1$ incident to $P$, and there is a 
line $\ell \in \mathcal C_2$ and a range $\mathcal R$ of points $P_1, \dots, P_k \in \mathcal C_2$ incident to $\ell$, such that the following 
condition is fulfilled:
\begin{enumerate}[(\rm CR)]
\item
There is a one-to-one correspondence between  $\mathcal P$ and $\mathcal R$ such that for any four lines $\ell_{i_1}, \ell_{i_2}, \ell_{i_3}, \ell_{i_4}
\in \mathcal P$, the cross ratio $(\ell_{i_1}, \ell_{i_2}; \ell_{i_3}, \ell_{i_4})$ is equal to the cross ratio $(P_{i_1}, P_{i_2}; P_{i_3},$ $P_{i_4})$ of 
the corresponding four points in $\mathcal R$.
\end{enumerate}
\end{defn}

Since duality is a one-to-one correspondence preserving cross-ratio, the following proposition is obvious. 

\begin{prop} \label{prop:polar_dual_comp}
Any configuration $\mathcal C$ is compatible with its dual $\mathcal C^{\delta}$.
\end{prop}

A simple consequence of Definitions~\ref{def:flex} and~\ref{def:comp} is as follows.

\begin{prop} \label{prop:FlexComp}
Two $k$-configurations are compatible if at least one of them is flexible. 
\end{prop}

\begin{proof}
Let $\mathcal C_1$ be a line-flexible configuration and $\ell \in \mathcal C_1$ a flexible line with a range of points $P_1, \dots, P_k \in \mathcal C_1$ 
incident to $\ell$. Let $\mathcal C_2$ be another configuration and choose a pencil of lines $\ell_1\dots, \ell_k \in \mathcal C_2$ incident to a point 
$P \in \mathcal C_2$. For all $i$, $i= 1, \dots, k$, take the points of intersection $M_i=  \ell \cap \ell_i$. Since $\ell$ is flexible, each $P_i$ can be 
chosen so as to coincide with $M_i$. Thus we have an \emph{elementary correspondence}~\cite{Cox1994} between the range of points $P_1, 
\dots, P_k$ and the pencil of lines $\ell_1\dots, \ell_k$, hence condition (CR) in Definition~\ref{def:comp} is fulfilled. By dualization, the same 
arguments apply in case $\mathcal C_1$ is point-flexible.
\end{proof}

Note that the converse is not necessarily true: there exist compatible configurations which are not flexible.


\subsection{Parametric affinity} \label{subsect:parametric}


Fix two different non-zero real numbers $a$ and $b$, and define%
\begin{equation} \label{eq:matrix}
\mathbf A_t =
\begin{matr}{cc}
\displaystyle\frac{a+t}{a} & 0 \\
0 & \displaystyle\frac{b+t}{b}
\end{matr},
\end{equation}

where $t \in \mathbb R$ is a parameter.

Fix a point $P_0$ with position vector $\mathbf r_0 = (x_0, y_0)$ lying on neither of the (Cartesian) coordinate axes. Multiplying by $\mathbf A_t$ 
transforms this point to a point %
\begin{equation*} 
\biggl(\frac{a+t}{a}x_0, \frac{b+t}{b}y_0\biggr) = (x_0, y_0) + t\biggl(\frac{x_0}{a}, \frac{y_0}{b}\biggr).
\end{equation*}
Thus, the set of all transforms of $P_0$ forms a straight line which passes through $P_0$ and whose 
direction vector is $(x_0/a, y_0/b)$. The parametric equation of this line is
\begin{equation} \label{eq:parametric}
(x,y) =  (x_0, y_0) + t\biggl(\frac{x_0}{a}, \frac{y_0}{b}\biggr),
\end{equation}

while the standard form of the equation of this line is
\begin{equation} \label{eq:intercept}
\frac{x}{\displaystyle\frac{a-b}{a}x_0} + \frac{y}{\displaystyle\frac{b-a}{b}y_0} = 1.
\end{equation}
Let $\alpha_t$ denote the affine transformation determined by the matrix $\mathbf A_t$; we shall use the family of affine transformations
\begin{equation} \label{eq:affinities}
\{\alpha_t \,|\, t \in \mathbb R\},
\end{equation}
which we denote by $\mathcal A(t)$.

\begin{defn}
For any point $P_0$ given by a position vector $\mathbf r_0(x_0, y_0)$, the line given by equation (\ref{eq:parametric}) (or equivalently, 
by (\ref{eq:intercept})) is called the \emph{orbit} of $P_0$ under the action of the family of transformations $\mathcal A(t)$.
\end{defn}

Here we use this term borrowed from the theory of group actions; we emphasize, however, that the analogy is only partial, since 
$\mathcal A(t)$ does not form a group.

By converse arguments, we obtain:

\begin{lem} \label{lem:orbit}
Let $L( \overline a, \overline b)$ be a straight line with intercepts $\overline a$ and $\overline b$ on the $x$- and $y$-axis, respectively, 
different from the origin. Then this line is the orbit of any of its points $P_0$ lying on none of the coordinate axes under the action of the 
parametric family of affine transformations given by a matrix of the form \rm{(\ref{eq:matrix})}, where
$a = -x_0 / \overline a$
and
$b= y_0 / \overline b$.
\end{lem}

\begin{proof}
The parametric equation of $L(\overline a, \overline b)$ can be written in the form 
$$
(x,y) = (x_0, y_0) + t((0, \overline b) - (\overline a, 0)) = (x_0, y_0) + t(-\overline a, \overline b).
$$
Comparing this equation with (\ref{eq:parametric}) verifies the assertion.
\end{proof}

Both this and the following lemma will be utilized in the next section.

\begin{lem} \label{lem:AffineProperties}
The affine transformation $\alpha_t$ has the following properties:

\begin{enumerate}[\rm (a)]
\item
$\alpha_t $ is the commuting product of two axial affinities determined by the matrices
\smallskip
\begin{equation} \label{eq:axial}
\begin{matr}{cc}
\displaystyle\frac{a+t}{a} & 0 \\
\phantom{\displaystyle\frac{t}{t}}0 \phantom{\displaystyle\frac{t}{t}} & 0
\end{matr}
\;\;\text{and}\;\;
\begin{matr}{cc}
0 & 0 \\
0 & \displaystyle\frac{b+t}{b}
\end{matr},
\end{equation}
\smallskip
whose axes are the $x$- and $y$-axis of the Cartesian coordinate system, respectively.
\item \label{intersection}
Given $k$ points $(x_0, y_i)$ ($i=0,1\dots, k$), lying on a line perpendicular to the $x$-axis, the corresponding orbits intersect in a common point on the $x$-axis. 
Similarly,  given $k$ points $(x_j, y_0)$ ($i=0,1\dots, k$), lying on a line perpendicular to the $y$-axis, the corresponding orbits intersect in a common point on 
the $y$-axis.
\end{enumerate}
\end{lem}
\begin{proof}
The commutation property is a trivial consequence of the fact that (\ref{eq:matrix}) is a diagonal matrix. The component matrices (\ref{eq:axial}) provide a well-known 
analytic description of axial affinities. Since $a$ and $b$ are fixed, the intercepts of line (\ref{eq:intercept}) depend only on $x_0$ and $y_0$; this verifies assertion
(\ref{intersection}).
\end{proof}

Comparing matrix $(2.1)$ in~\cite{BerGevPis2021} with our matrix (\ref{eq:matrix}) given here shows that the latter is 
a generalization of the former. In fact, choosing $a=-b=h$ returns matrix $(2.1)$ with the integer values $t = j = 1, \dots , h-1$.

\begin{figure}[!h]
\begin{center}
\includegraphics[width=.7\textwidth]{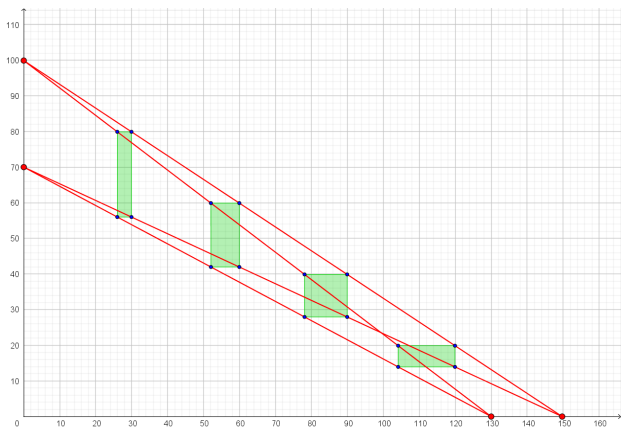}
\caption{Illustration of the action of affine transformations $\alpha_t$ on a square at four distinct values of $t$. The original square is 
              $t=0$, and the affine transformations applied use $t = -2, -1, 1$). Orbits of vertices of the square are indicated by red segments.}
\label{fig:rectangles}
\end{center}
\end{figure}

Figure~\ref{fig:rectangles} shows an example of the action of $\alpha_t$ on a square with side length 12. Here $a=-1.5b$, and the upper right vertex of 
the square is the point $(90, 40)$. The square is obviously fixed at $t=0$. Choosing $a=2$, the non-identical transforms of the square shown in the 
figure have values $t = -2, -1, 1$. The points of intersection of suitable pairs of vertex orbits are highlighted by red; they lie in pairs on the 
coordinate axes, illustrating Lemma~\ref{lem:AffineProperties}(\ref{intersection}).


\section{Gr\"unbaum incidence calculus extended} \label{sect:IncCalc}


In \cite{Gru2009b}, Section 3.3, Branko Gr\"unbaum described a sequence of operations to produce new $(n_4)$ configurations from various input configurations. 
These operations were called the ``Gr\"unbaum Incidence Calculus'' in \cite{PisSer2013}, Section 6.5. Some of the operations Gr\"unbaum described are specific
to producing 3- or 4-configurations. Other operations can be generalized in a straightforward way to produce $(n_k)$ configurations either from  smaller $(m_k)$
configurations with certain properties, or from $(m_{k-1})$ configurations.

In \cite{BerGevPis2021}, we generalized two constructions, the \emph{affine replication} and \emph{affine switch} operations. In this section we revise them, the latter in 
more detail, generalizing it by using parametric affinities. We also generalize three other constructions, the \emph{parallel switch} and two \emph{deleted union} operations. 
In the next two sections, we use these constructions to produce a large number of relatively small 5- and 6-configurations, in order to decrease the bounds on 
$N_{5}$ and $N_{6}$. 


\subsection{Affine replication} \label{subsect:AR}

This operation, generalizing Gr\"unbaum's $(\mathbf{5m})$ construction, takes as input an $(m_{k-1})$ configuration $\mathcal C$ and produces a 
$((k+1)m_k)$ configuration $\mathcal D$ with a pencil of $m$ parallel lines. A sketch of the construction is that we apply to $\mathcal C$ a sequence 
of $k$ orthogonal axial affinities $\alpha_1, \dots, \alpha_k$ with a suitably chosen common axis $A$. Then, each point $P$ of $\mathcal C$ and its 
$k$ images will be collinear. Moreover, each line $\ell$ of $\mathcal C$ and its $k$ images are concurrent at a single point, a fixed point of the affinities, 
lying on the common axis. The new configuration $\mathcal D$ consists of the points and lines of $\mathcal C$ and its images, the new lines corresponding 
to the collinearities from each point $P$, and the new points corresponding to the concurrences from each line $\ell$. A detailed description of this construction 
is given in~\cite{BerGevPis2021} (see also the examples given there in Figures 2 and 3). It is denoted by $\AR(m_k)$. 

\begin{prop} \label{prop:AR}
If affine replication $\AR(m_k)$ is applied to any $(m_{k-1})$ configuration $\mathcal C$, the result is a $(((k + 1)m)_k)$ line-flexible configuration 
with a pencil of $m$ parallel flexible lines. The flexible lines are precisely the lines connecting a configuration point with each of its affine images.
\end{prop}
\begin{proof}
Let $P$ be any point of $\mathcal C$ at a distance $d$ from the axis $A$. Let $r_1, \dots, r_k$ be the ratios of the affinities $\alpha_1, \dots, \alpha_k$. 
Then the distance of the points $\alpha_1(P), \dots, \alpha_k(P)$ from $A$ is $r_1d, \dots, r_kd$, respectively; moreover, all these points lie on a common 
line $\ell_P$ perpendicular to $A$ (recall the basic properties of an orthogonal axial affinity). Hence it is clear that using suitable affinities, the position of the 
points $P, \alpha_1(P), \dots, \alpha_k(P)$ on the line $\ell_P$ can be prescribed. Thus the new configuration is line-flexible, and the flexible lines are 
the $m$ parallel new lines. The rest of the proposition is contained precisely in Lemma 2.2 in \cite{BerGevPis2021}.
\end{proof}


\subsection{Affine switch} \label{subsect:AS}

The original version of this operation, denoted by $(\mathbf{3m+})$ and applied to $(m_4)$ configurations, occurred also in Gr\"unbaum's 
book \cite{Gru2009b}, Section 3.3. In~\cite{BerGevPis2021}, it was generalized to $(m_k)$ configurations, but the affine transformations 
used in that presentation were parameterized by integers. Here we generalize it further, using a continuous parameter; this generalization results 
in line-flexible configurations.

Suppose that $\mathcal C$ is an $(m_k)$ configuration that contains a pencil $\mathcal P$ of parallel lines $\ell^1, \dots, \ell^p$ with $p\ge 1$ and a pencil 
$\mathcal Q$ of parallel lines $\ell^{p+1}, \dots,\ell^{p+q}$ with $q \ge 0$ such that the two pencils are \emph{independent}, which means that they share 
no common configuration points. We also assume that these pencils are perpendicular to each other, since if the configuration has non-perpendicular pencils, it 
can easily be transformed, using a suitable affine transformation, into to an isomorphic copy in which two arbitrarily chosen independent parallel pencils will be 
perpendicular. 

We slightly change the notation in~\cite{BerGevPis2021} and use $\AS(m_{k}, r_p, r_q)$ to denote this construction, where  $r_p$ and $r_q$ denote the number 
of deleted lines from pencil $\mathcal P$ and $\mathcal Q$, respectively. 

\begin{prop}
Starting from any $(m_k)$ configuration with independent pencils of $p \ge 1$ and $q \ge 0 $ parallel lines, for each integer $r$ with $1 \le  r \le p+q$, the 
$\AS(m_{k}, r_p, r_q)$ construction produces a line-flexible $(n_k)$ configuration, where $n = (k-1)m + r$. This configuration contains independent pencils 
with $p'=(k-1)(p-r_p)$ and $q'=(k-1)(q-r_q)$ parallel lines, respectively.
\end{prop}

\begin{proof}

The steps of the construction are as follows.
\begin{enumerate} 
\item
Take a copy of $\mathcal C$ in a position such that $\mathcal P$ and $\mathcal Q$ are parallel with the $x$- and $y$-axis of the Cartesian 
coordinate system, respectively, and none of the configuration points is incident to either of the coordinate axes.
\item
Take the line $\ell^1$ of $\mathcal P$, choose a configuration point $P_1$ that is incident to this line, and draw a straight line $L(P_1)$ 
through $P_1$ such that it intersects (say) the positive branches of the coordinate axes, avoids all the other configuration points of 
$\mathcal C$, and is not parallel with any line connecting two points of $\mathcal C$. 
\item
Choose $k-2$ points (arbitrarily) on $L(P_1)$ different from $P_1$ and from each other such that none of them lie on either of the coordinate axes.

By Lemma~\ref{lem:orbit}, each of these $k-2$ points can be considered as the transforms of $P_1$ by affinities $\alpha_t$ 
with $k-2$ different values of $t$. By a slight abuse of notation, we shall denote these affinities $\alpha_{t_1}, \dots, \alpha_{t_{k-2}}$ 
by $\alpha_{1}, \dots, \alpha_{k-2}$, respectively. Thus we have the $k-1$ points 
\begin{equation*} \label{eq:PointTansforms}
P_1, \alpha_1(P_1), \dots, \alpha_{k-2}(P_1).
\end{equation*}
\item
Using the pairs $(P_1, \alpha_1(P_1)), \dots, (P_1, \alpha_{k-2}(P_1))$, construct the images 
\begin{equation} \label{eq:images}
\mathcal C_1= \alpha_1(\mathcal C), \dots, 
\mathcal C_{k-2}= \alpha_{k-2}(\mathcal C).
\end{equation}

Note that the affinities $\alpha_t$ are products of axial affinities (cf.\ Lemma~\ref{lem:AffineProperties}(a)); moreover, an axial 
affinity is determined by its axis and a pair consisting of a point and its image. Thus the constructions above can easily be realized 
(recall e.g.\ Proposition 2.8 in~\cite{BerGevPis2021} for the details). 
\item
Take now the configuration point $P_2$ incident to $\ell^1$, and consider its images 
\begin{equation*} \label{eq:PointTansforms}
\alpha_1(P_2), \dots, \alpha_{k-2}(P_2).
\end{equation*}
Again by Lemma~\ref{lem:orbit}, these points will lie on the orbit of $P_2$, which we denote by $L(P_2)$.
\item
Repeat the preceding step for the rest of the configuration points $P_3, \dots, P_k$ incident to $\ell^1$, and denote the corresponding lines 
by $L(P_3), \dots, L(P_k)$.  By Lemma~\ref{lem:AffineProperties}(b), these lines, along with $L(P_{1})$ and $L(P_{2})$, intersect in a 
common point on the $y$-axis. Denote this point by $Y^1$.
\item
Remove $\ell^1$ from $\mathcal C$ and all its images from the images of $\mathcal C$ given in~\eqref{eq:images}. We denote 
the resulting structure by  $\mathcal C'$.
\item
Form a new configuration as the following union:
$$
\{\mathcal C',  \alpha_1(\mathcal C'), \dots, \alpha_{k-2}(\mathcal C')\} \cup \{L(P_1), \dots, L(P_k)\} \cup \{Y^1\}.
$$
\end{enumerate}

The configuration obtained in this way is a $(((k-1)m+1)_k)$ configuration.

The procedure can be repeated at most $p-1$ times using further lines $\ell^2, \dots, \ell^p$ of $\mathcal P$.
Furthermore, interchanging the role of $\mathcal P$ and $\mathcal Q$ as well as of the coordinate axes, the 
whole procedure can be repeated again, with the same affine images as before, but using the lines of $\mathcal Q$
and the configuration points incident to them. In this way we obtain at most $r=p+q$ distinct new $k$-configurations.

Observe that with reference to Lemmas~\ref{lem:orbit} and~\ref{lem:AffineProperties} where it is needed,
realizability of all the 8 steps above is verified. In addition, each of the $r$ new configurations that can obtained 
in this way is line-flexible, which is a simple consequence of the freedom in choosing the $k-1$ points in steps 
(2) and (3). 

The rest of the proposition on independent pencils of parallel lines is straightforward.
\end{proof}
\begin{figure}[!h]
\begin{center}   
        \includegraphics[width=0.95\textwidth]{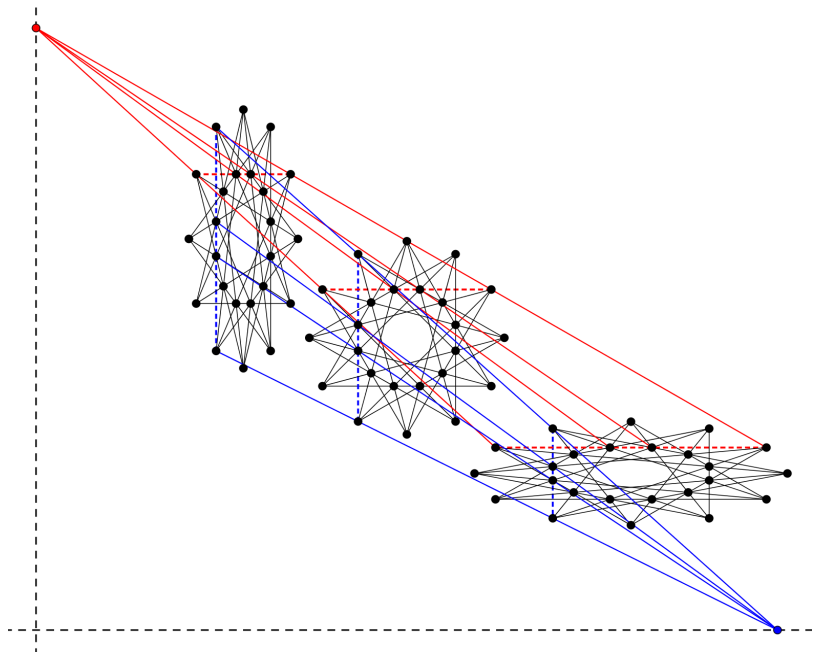}
\end{center} 
\caption{A $(74_{4})$ configuration formed by applying the affine switch operation to a $(24_{4})$ configuration. The $(24_{4})$ configuration has two independent pencils, but only one line is used from each pencils, for clarity. The dashed lines are excluded from the final configuration.}
\label{fig:AS-Example}
\end{figure}

Figure~\ref{fig:AS-Example} shows example of the affine switch operation applied to a configuration of type $(24_{4})$. In the construction, the configuration has two independent pencils which are perpendicular to each other and which both consist of two lines. In the figure, only one line is used from each of the pencils, producing a $(74_4)$ configuration.
(For an example where consecutively one, two and three lines of a single pencil are used, see Figure 5 in~\cite{BerGevPis2021}.)


\subsection{Parallel switch} \label{subsect:PS}


Start from an $(m_k)$ configuration $\mathcal C$, and delete a line from it. We denote the incidence structure obtained in this way by $\overline{\mathcal C}_0$.
Take $k-1$ translated copies $\overline{\mathcal C}_1,\dots, \overline{\mathcal C}_{k-1}$, where the direction of translation differs from that of each line of 
$\mathcal C$; the $k-1$ distances of translation are arbitrary, but must be chosen so that no new incidences occur. Finally, add $k$ new parallel lines 
to the set-union $\overline{\mathcal C}_0 \cup \overline{\mathcal C}_1 \cup \dots \cup \overline{\mathcal C}_{k-1}$ such that they pass through the 
``defective" points of $\overline{\mathcal C}_0$ and of its $k-1$ translates. We denote this construction as $\PS(m_{k})$ and note that this construction 
also occurs in~\cite{PisSer2013}. 

\begin{figure}[!h]
\begin{center}
\includegraphics[width = .225\textwidth]{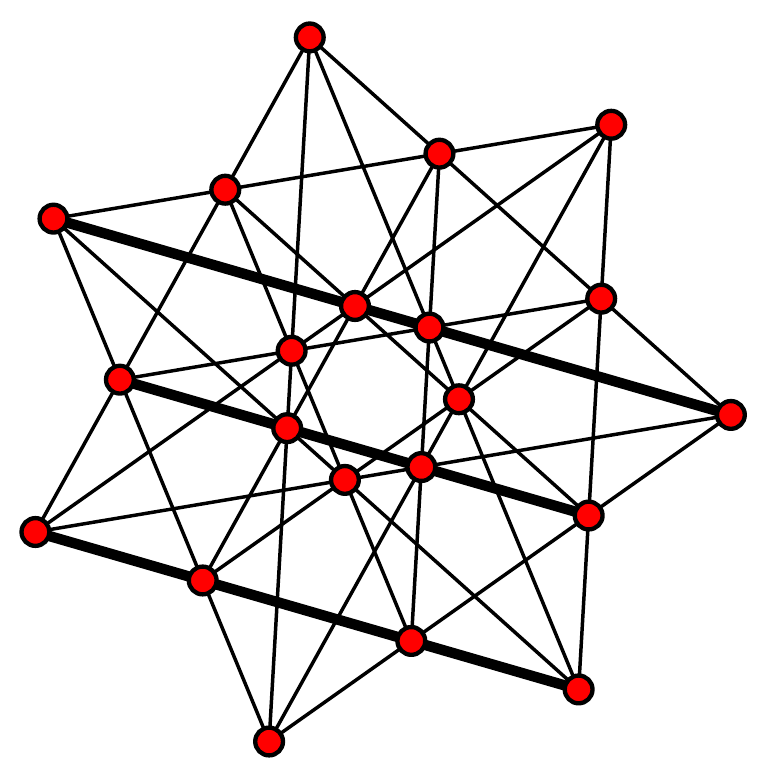}
\includegraphics[width = \textwidth]{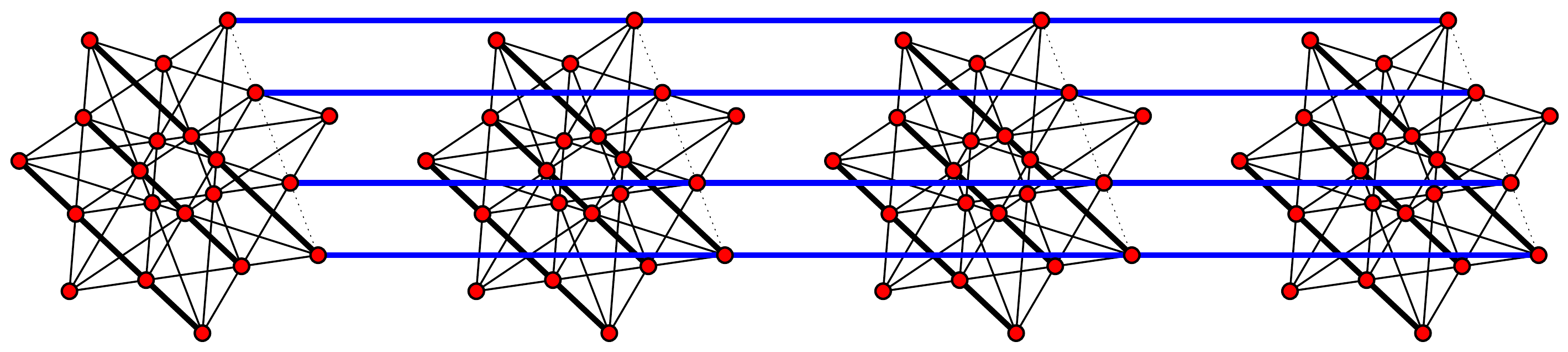}
\caption{The $(21_{4})$ Gr\"unbaum--Rigby configuration $\mc{C}$ (top), and the parallel switch construction applied to it (bottom). The lightly dotted lines 
have been removed. The construction produces a $(84_{4})$ configuration, which has a parallel pencil of 4 new lines generated by the construction. However, 
since $\overline{\mc{C}}$ contained a pencil of 3 parallel lines already, $\mc{D}$ also contains a pencil of 12 parallel lines (black). }
\label{fig:ParallelSwitch}
\end{center}
\end{figure}

As a result of the $\PS(m_{k})$ construction, a $((km)_k)$ configuration is obtained which contains a parallel pencil of $k$ new lines. However, any parallel pencil 
in $\overline{\mathcal C}_0$ is multiplied by $k$ (see the example given in Figure~\ref{fig:ParallelSwitch}) as well, which produces at least as many parallel lines 
as the new ones, and usually more.
\begin{prop}
Start from an $(m_k)$ configuration $\mathcal C$. Assume that it has a pencil $\mathcal P$ of $p$ lines, and an independent pencil  $\mathcal Q$ of $q$ lines, 
none of which is the line deleted in the construction process. Then the parallel switch construction $\PS (m_k)$ produces from $\mathcal C$ a line-flexible $((km)_k)$ 
configuration $\mathcal D$. This configuration contains a pencil of $kp$ lines, and an independent pencil of $kq$ lines. 
\end{prop}


\subsection{The deleted union constructions} \label{subsect:du}


These constructions again generalize a construction of Gr\"unbaum, which he briefly described in~\cite[pp.\ 180--182]{Gru2009b}.


\subsubsection{The $\du(\mathcal C_{1}, \mathcal C_{2})$ construction} \label{subsubsect:basic}


This is the basic case of the deleted union constructions, and each of the other cases can be considered as a modification of this. 
Here $\mathcal C_{1}$ and $\mathcal C_{2}$ are two compatible $k$-configurations, usually different  (cf.\ Definition \ref{def:comp}). 
Let  $\mathcal P$ be a pencil of lines in $\mathcal C_{1}$ incident to a point $P \in \mathcal C_{1}$, and let $\mathcal R$ be a 
range of points in $\mathcal C_{2}$ incident to a line $\ell \in \mathcal C_{2}$, defining the compatibility of these configurations. 
Apply a suitable projective collineation $\pi$ to $\mathcal C_{2}$ which produces a configuration $\mathcal C'_{2}$ such that each 
point in $\mathcal R' \subset \mathcal C'_{2}$ is incident to a corresponding line in $\mathcal P$, but is not incident to any other 
line in $\mathcal C_{1}$. Finally, delete the point incident to the lines of $\mathcal P$, and delete the line incident to the points of 
$\mathcal R'$. 

\begin{figure}[!h]
\begin{center}
\includegraphics[width = .33 \linewidth]{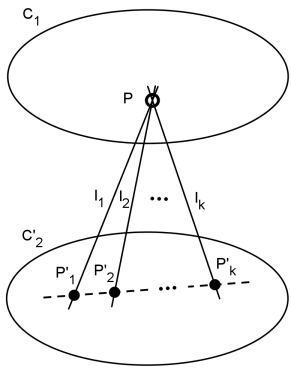}
\caption{A scheme illustrating the $\du(\mathcal C_{1}, \mathcal C_{2})$ construction.  
             The point and line deleted in the last step are denoted by an empty circle and a dashed line, respectively.}
\label{fig:DUscheme}
\end{center}
\end{figure}

This construction is illustrated by the scheme in Figure~\ref{fig:DUscheme}.

The following proposition is straightforward.

\begin{prop} \label{prop:du(c1c2)}
Let $\mathcal C_1$ and $\mathcal C_2$ be two compatible configurations with the respective symbols $(n_k)$ and $(n'_k)$. 
Then $\du(\mathcal C_1, \mathcal C_2)$ exists and is a configuration with symbol $((n+n'-1)_k)$.
\end{prop}

Note that a suitable projective collineation could equally well be applied to $\mathcal C_{1}$ instead of $\mathcal C_{2}$, thus yielding an isomorphic 
copy of the configuration obtained in the former case.

The crucial part the $\du(\mathcal C_{1}, \mathcal C_{2})$ construction (and in each of its modified versions) is to find the precise form of the projective 
collineation $\pi$. We proceed in the following steps.
\begin{enumerate}
\addtocounter{enumi}{-1}
\item
Let $\mathcal P = \{\ell_1, \ell_2,\dots,\ell_k\} \subset \mathcal C_{1}$ be a pencil of lines used in the construction (cf.\ Figure~\ref{fig:DUscheme}).
Let $\ell'$ be a line incident to none of the points of $\mathcal C_{1}$ or $\mathcal C_{2}$. Take the intersections $P'_i = \ell' \cap \ell_i$ for all
$i= 1, \dots , k$.
\item
Let $\mathcal R = \{P_1, P_2,\dots,P_k\} \in \mathcal C_{2}$ be a range of points used in the construction. The first main step of our procedure
is to establish a one-dimensional projectivity $\pi^{(1)}$ between this range and the range $\{P'_1, P'_2, \dots, P'_k\} \subset \ell'$ as follows:
\begin{equation}
P_1 P_2 \cdots P_k \barwedge P'_1 P'_2 \cdots P'_k .
\end{equation}
This can be realized by the product of two perspectivities:
\begin{equation} \label{eq:twopersp}
P_1 P_2 \cdots P_k \doublebarwedge Q_1 Q_2 \cdots Q_k \doublebarwedge P'_1 P'_2 \cdots P'_k,
\end{equation}
whose centres are (say) $P'_1$ and $P_1$, respectively (cf.\ Figure 1.6B and relationship (1.61) in~\cite{Cox1994}). (Here we use the traditional 
transformation symbols due to von Staudt and Veblen~\cite{Cox1994}.)
\item
In this second main step the projectivity $\pi^{(1)}$ is extended to a two-dimensional projectivity $\pi := \pi^{(2)}$, which is a projective collineation
transforming the whole configuration $\mathcal C_{2}$ to $\mathcal C'_{2}$:
\begin{equation}
\pi: \mathcal C_{2} \mapsto \mathcal C'_{2}.
\end{equation}
This can be realized by extending each of the two perspectivities in (\ref{eq:twopersp}) to a perspective collineation. The centre of the first and second 
collineation is the point $P'_1$ and $P_1$, respectively. When chosing the axis for these collineations, observe that there is a point $F'=\ell' \cap q$ different 
from $P'$ fixed by the first perspectivity in (\ref{eq:twopersp}), and likewise a point $F=\ell \cap q$ different from $P$ fixed by the second perspectivity; 
hence the axis of the first collineation must go through $F'$, while the axis of the second collineation must go through $F$.

\end{enumerate}

We note that the (technically, somewhat lengthy) procedure described above can be omitted, and thus the {\bf DU} construction can be made 
much simpler, in the following special case. 

\begin{rem} \label{rem:flex}
Assume that one of the component configurations is line-flexible. In this case we apply Proposition~\ref{prop:FlexComp}, 
and the precise position of the points incident to the line to be deleted can easily be found using the proof of this proposition.
(Note that the dual case, i.e.\ if one of the component configurations is point-flexible, can be treated similarly.)
\end{rem}

As a consequence, we have the following special case of the previous proposition. 

\begin{prop} \label{prop:flex}
Let $\mathcal C_1$ and $\mathcal C_2$ be two configurations with the respective symbols $(n_k)$ and $(n'_k)$ and let $\mathcal C_2$ be flexible. 
Then $\du(\mathcal C_1, \mathcal C_2)$ exists and is a configuration with symbol $((n+n'-1)_k)$.
\end{prop}

\begin{figure}[!h]
\begin{center}
\includegraphics[width = .66\linewidth]{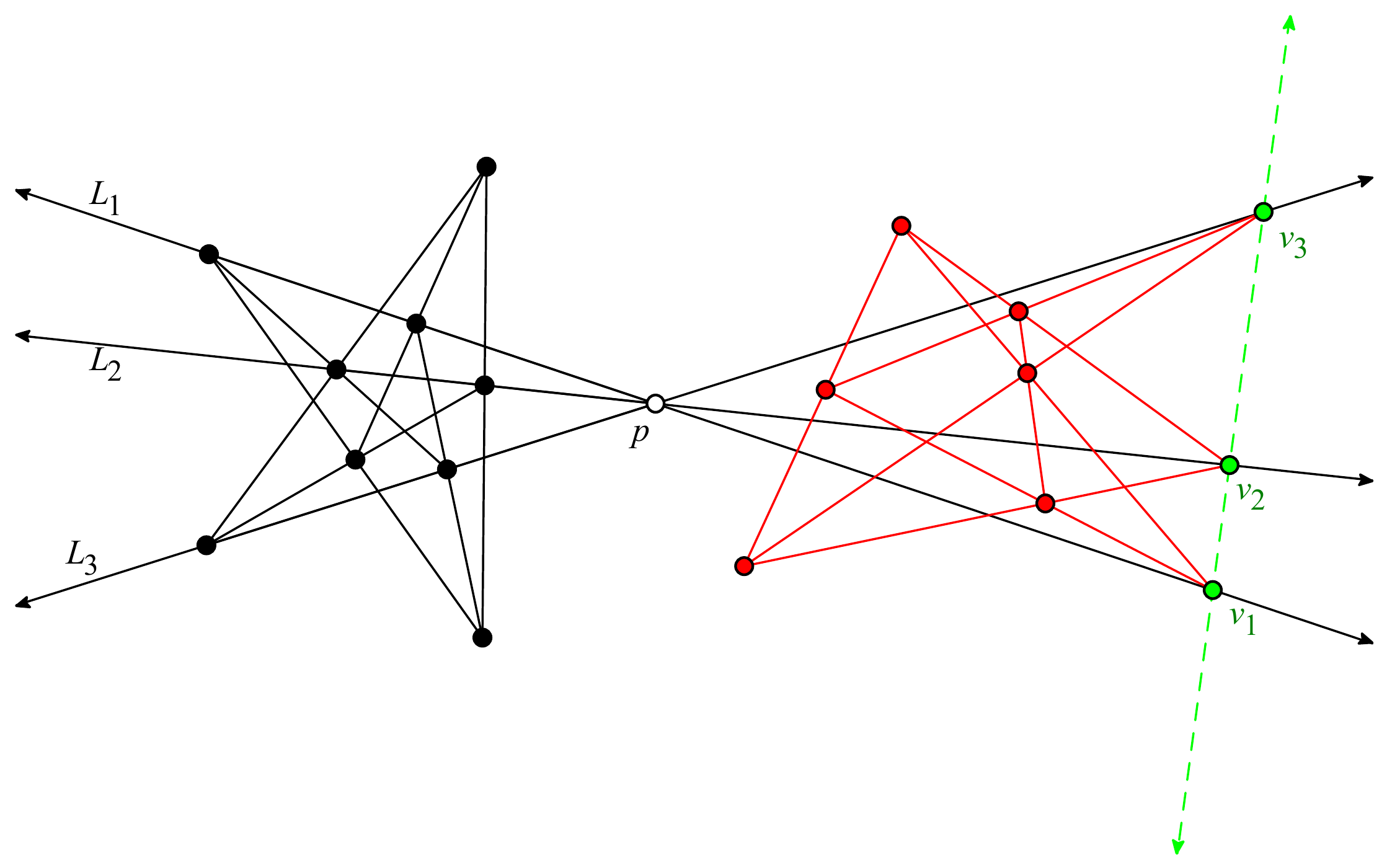}
\caption{Example of the $\du(\mathcal C_{1}, \mathcal C_{2})$ construction, where $\mathcal C_{1}$ is the cyclic $(10_{3})$ configuration 
              and $\mathcal C_{2}$ is the $(9_3)$ Pappus configuration. The construction produces an $(18_{3})$ configuration. Neither the white point nor the green dashed line are elements of the final configuration.}
\label{fig:DU(10-3-9-3)}
\end{center}
\end{figure}

Figure~\ref{fig:DU(10-3-9-3)} shows an example of this special case of the construction, applied to a $(10_3)$ and a $(9_3)$ configuration. 
Here the second component, the Pappus configuration, is a flexible configuration; its flexible line used in the construction is highlighted by 
green. The white point $p$ and the dashed green line is deleted in the last step of the construction.


\subsubsection{The $\du(1)(\mc{C})$ construction} \label{sect:du(1)}


When discussing the deleted union construction in case of $(n_4)$ configurations, Gr\"{u}nbaum observes that one can choose for $\mathcal C_2$ the polar 
of $\mathcal C_1$~\cite{Gru2009b}. This observation is naturally valid for any $(n_k)$ configuration in a slightly more general form, thus we introduce for 
the $\du(\mathcal C, \mathcal C^{\delta})$ construction the notation $\du(1)(\mc{C})$, or simply $\du(1)$ if we do not want to specify the configuration. 
Indeed,

We have the following simple consequence of Propositions  \ref{prop:polar_dual_comp} and \ref{prop:du(c1c2)}.

\begin{cor} \label{cor:du(1)}
Let $\mathcal C$ be any $(n_k)$ configuration. Then  
$\du(1)(\mc{C})$ exists and is a configuration with symbol $((2n-1)_k)$.
\end{cor}

The polarity applied here is induced by a conic~\cite{Cox1994}. Often, when $\mc C$ is a configuration with rotational symmetry, we take this conic 
to be a circle concentric with the centre of rotation; hence the polarity in this particular case reduces to a reciprocation with respect to the circle in 
question (cf.\ e.g.~\cite{CoxGre1967}). 

\begin{figure}[!h]
\begin{center}
\includegraphics[width=.6\linewidth]{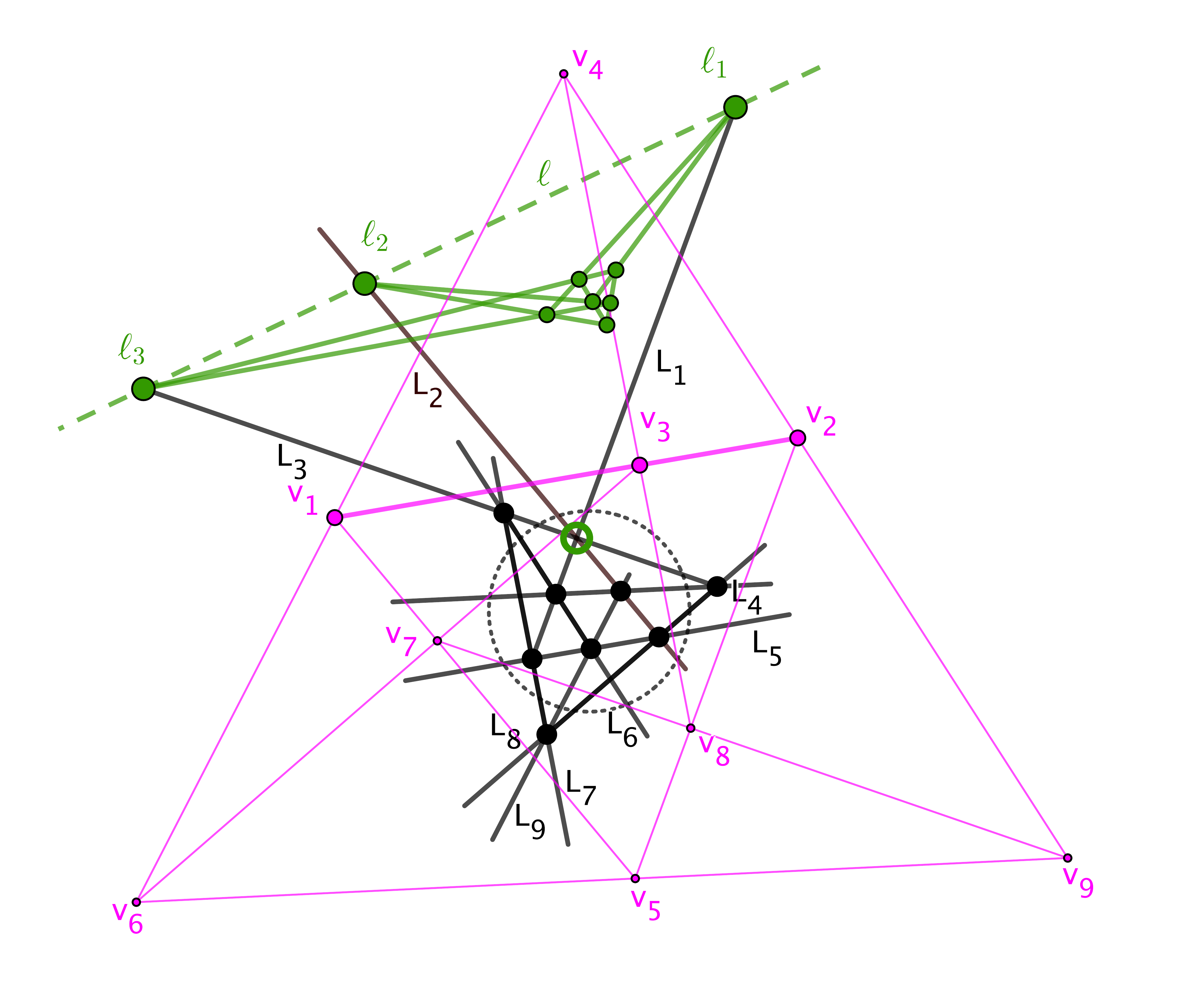}
\caption{The $\du(1)$ construction applied to one of the non-Pappus $(9_{3})$ configurations (for the details, see the text). 
}
\label{fig:9-3-DU1-details}
\end{center}
\end{figure}

An example of this construction is shown in Figure~\ref{fig:9-3-DU1-details}. The original configuration $\mc C$ is shown in black, with the deleted point 
$p$ shown as a hollow green circle and the lines of the configuration labelled $L_{i}$. The dotted circle is the circle of reciprocation, and the magenta points 
$v_{i}$ are the poles of the lines $L_{i} \in \mc C$ obtained under the reciprocation. The arbitrary line $\ell$ intersecting the lines $L_{1}$, $L_{2}$, $L_{3}$ 
is shown dashed green. The projective collineation $\pi$ mapped $v_{1}, v_{2}, v_{3}$ in order to $\ell_{1}, \ell_{2}, \ell_{3}$, and the rest of the points 
$v_{i}$ were mapped using this collineation to the other green points (unlabeled for clarity). The final configuration $\du(1)(\mc C)$ is the union of the solid 
black and green points and lines.

\subsubsection{The $\du(t)$ construction} \label{subsubsect:du(t)}

This iterates the $\du(1)$ construction, except that instead of working with a single  configuration point $p$ and arbitrary line $\ell$, 
the construction uses several points $p^{1}, p^{2}, \ldots, p^{t-1}$ and lines $\ell^{1}$, $\ell^{2}, \ldots, \ell^{t-1}$. 
That is, we construct each of the configurations $(\mc{C}^{i})'$ associated with each pair $p^{i}$, $\ell^{i}$ 
independently, and then delete all the points $p^{i}$ and lines $\pi_{i}(\varrho(p^{i}))$ and combine what 
remains (here $\varrho$ denotes the polarity, which in some cases reduces to reciprocity). 

It can directly be seen that for the iteration any number of configuration points can be used (each at one time). Thus we have:

\begin{prop}
Let $\mc C$ be any $(n_k)$ configuration and $t \le n$ a positive integer. Then $\du(t)(\mc C)$ exists and is a configuration with symbol $(((t+1)n-t)_k)$. 
\end{prop}

\begin{figure}[!h] 
\begin{center} 
   \includegraphics[width=.55\linewidth]{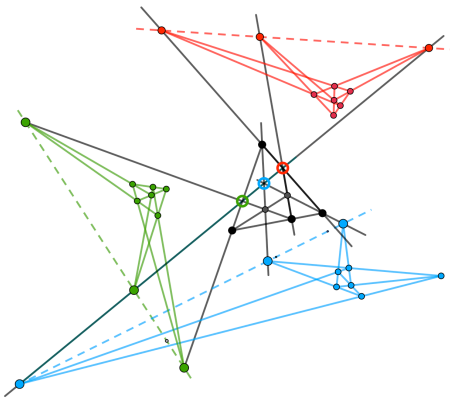}
   \caption{The configuration $\du(3)(9_{3})$ is a $(33_{3})$ configuration. 
                 Hollow points have been deleted from the original configuration (shown in black in the middle), and the dashed lines are 
                 the lines $\ell^{1}$, $\ell^{2}$, $\ell^{3}$ which are not lines of the final configuration.}
   \label{fig:9-3-DU3}
\end{center}
\end{figure}

Figure \ref{fig:9-3-DU3} shows a $(33_{3})$ configuration constructed by applying the $\du(3)$ construction to a $(9_{3})$ configuration: the construction 
can be applied independently multiple times, even to points that are collinear in the original construction. (We note that the non-Pappus $(9_3)$ configuration 
used here and in the previous figure is the one specified as $(9_3)_2$ by Hilbert and Cohn-Vossen~\cite{HCV}.)

\begin{rem} \label{rem:pencils}
Since projective transformations do not preserve parallel lines, no additional pencils of parallel lines are guaranteed to be generated by this construction. 
Thus, if a configuration $\mc{C}$ has independent pencils with $p$ and $q$ parallel lines, the configuration $\mc{\hat C}$ arising from applying $\du(t)$ 
to $\mc{C}$ will also have independent pencils of $p$ and $q$ parallel lines.
\end{rem}


\subsubsection{The $\du^{(2)}(\mc C_{1}, \AR(\mc C_{2}))$ construction}\label{sec:DU-2}


This is a modified version of the $\du(\mc C_{1}, \mc C_{2})$ construction;  here we use two points of $\mathcal C_1$ and two lines of 
$\AR(\mathcal C_{2})$ to be deleted in the last step of the construction (the superscript in the notation refers to this feature). (We note that such
a modification of the $\du$ construction is also used by Gr\"unbaum, see \cite[Figure 3.3.18]{Gru2009b}.) 

Assume that $\mc C_1$ is a $k$-configuration. Then Proposition~\ref{prop:AR} implies that $\mc C_2$ must be a $(k-1)$-configuration; moreover,
$\AR(\mc C_{2})$ is a line-flexible configuration. Observe that if we used a simple $\du$ construction (with one deleted point and one deleted line),
then (by Remark~\ref{rem:flex}) the latter flexibility property would imply that the construction could easily be realized without any further condition.
However, in the present case some additional conditions are needed, as follows.

{\sc Condition$\!$} (P). Configuration $\mathcal C_1$ must contain a pencil $\mc P$ through a configuration point $D$ and a pencil $\mc P'$ 
through a configuration point $D'$ such that these two pencils form a centrally symmetric pair with respect to a centre $O$ (which is not necessarily a 
configuration point).

The second condition concerns the mutual position of $\mc C_{1}$ and $\AR(\mc C_{2})$ guaranteeing that the construction is feasible.

{\sc Condition$\!$} (M). Let $P_0$ and $P'_0$ be two configuration points of $\mc C_{2}$, and let $M$ be the midpoint of the segment $P_0P'_0$. 
Then choose a mutual position of $\mc C_{1}$ and $\mc C_{2}$ such that 

\hskip 15 pt (M1) $M$ coincides with $O$;

\hskip 15 pt (M2) the axis of affinity used in the $\AR(\mc C_{2})$ construction goes through $M$;

\hskip 15 pt (M3) one of the lines in pencil $\mc P$ must go through $P_0$ (denote it by $\ell_0$).

\begin{figure}[!h]
\begin{center}
\includegraphics[width=.66\textwidth]{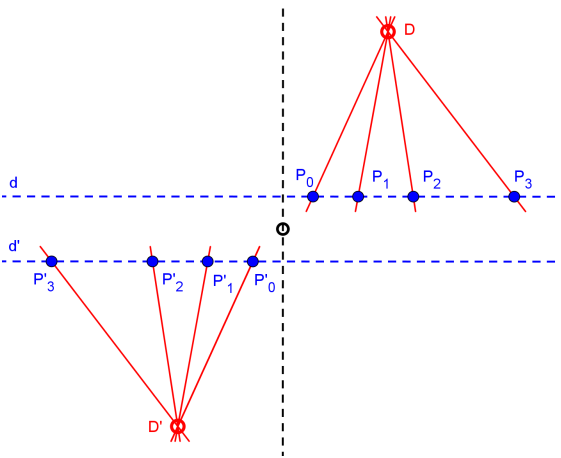}
\caption{Illustration for the $\du^{(2)}(\mc C_{1}, \AR(\mc C_{2}))$ construction with $k=4$. Elements belonging to $\mc C_{1}$ are coloured red, while those 
belonging to $\AR(\mc C_{2})$ are distingushed by blue. The vertical black line is the axis of affinity, and the centre of symmetry is denoted by an empty black circle. 
}
\label{fig:DU_AR}
\end{center}
\end{figure}

Now observe that the symmetry conditions ensure that the symmetric counterpart of the line $\ell_0$ will go through $P'_0$. Furthermore, recall that the flexible lines of 
$\AR(\mc C_{2})$ are precisely the lines connecting a configuration point with each of its affine images (cf.\ Proposition~\ref{prop:AR}). In particular, the line connecting 
the points $P_0, P_1=\alpha_1(P_0), P_2=\alpha_2(P_0), \dots, P_k=\alpha_k(P_0)$, where $\alpha_i$ denotes the $i$th affine transformation, will be a flexible line 
of $\AR(\mathcal C_{2})$. For all $i$ $(i= 1,2, \dots, k)$, the transformation $\alpha_i$ must be chosen so that $\alpha_i(P_0)$ is incident with the $i$th line of pencil 
$\mc P$. Practically, this means that first we determine the points $P_i$ as intersections $\ell_i \cap d$ where $\ell_i$ is the $i$th line in the pencil $\mc P$ and $d$ is 
a line going through $p_0$ and perpendicular to the axis of affinity. Then, by this intersection condition, all the $\alpha_i$ affinities are determined (recall that an axial 
affinity is determined by its axis and one pair of corresponding points; hence in our case $\alpha_i$ is determined by the pair $(P_0, P_i)$ for all $i$). For an illustration,
see Figure~\ref{fig:DU_AR}.

The symmetry conditions imposed above ensure that for the points $P_i$ and the line $d$ we have their symmetric counterparts (with respect to the centre $O$) $P'_i$ 
and $d'$ such that the points are incident to the corresponding lines of $\mc P'$, and $d'$ is a (flexible) line of $\AR(\mc C_{2})$ parallel to $d$.  Thus, what remains 
is to delete the points $Q$, $Q'$ and the lines $d$, $d'$, and the construction $\du^{(2)}(\mc C_{1}, \AR(\mc C_{2}))$ is now complete.

\begin{figure}[!ht]
\begin{center}
\ffigbox[]{
\begin{subfloatrow}[3]
\ffigbox{\caption{$\mc C_{1}$.}\label{fig:DU2start}}{
    \includegraphics[width=.3\textwidth]{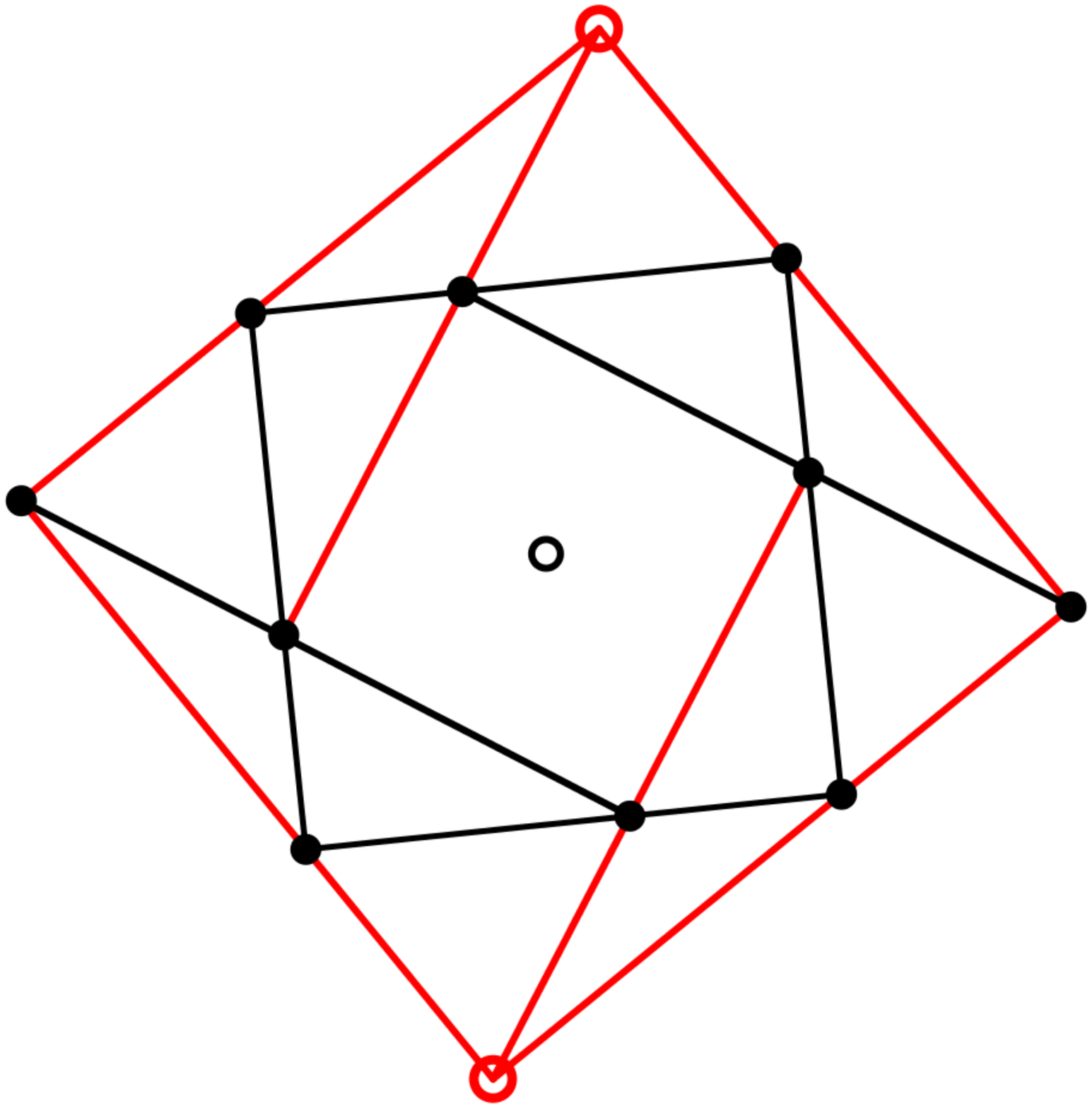}
    } 
    \ffigbox{\caption{$\AR(\mc C_{2})$.}\label{fig:DU2middle}}{
    \includegraphics[width=.3\textwidth]{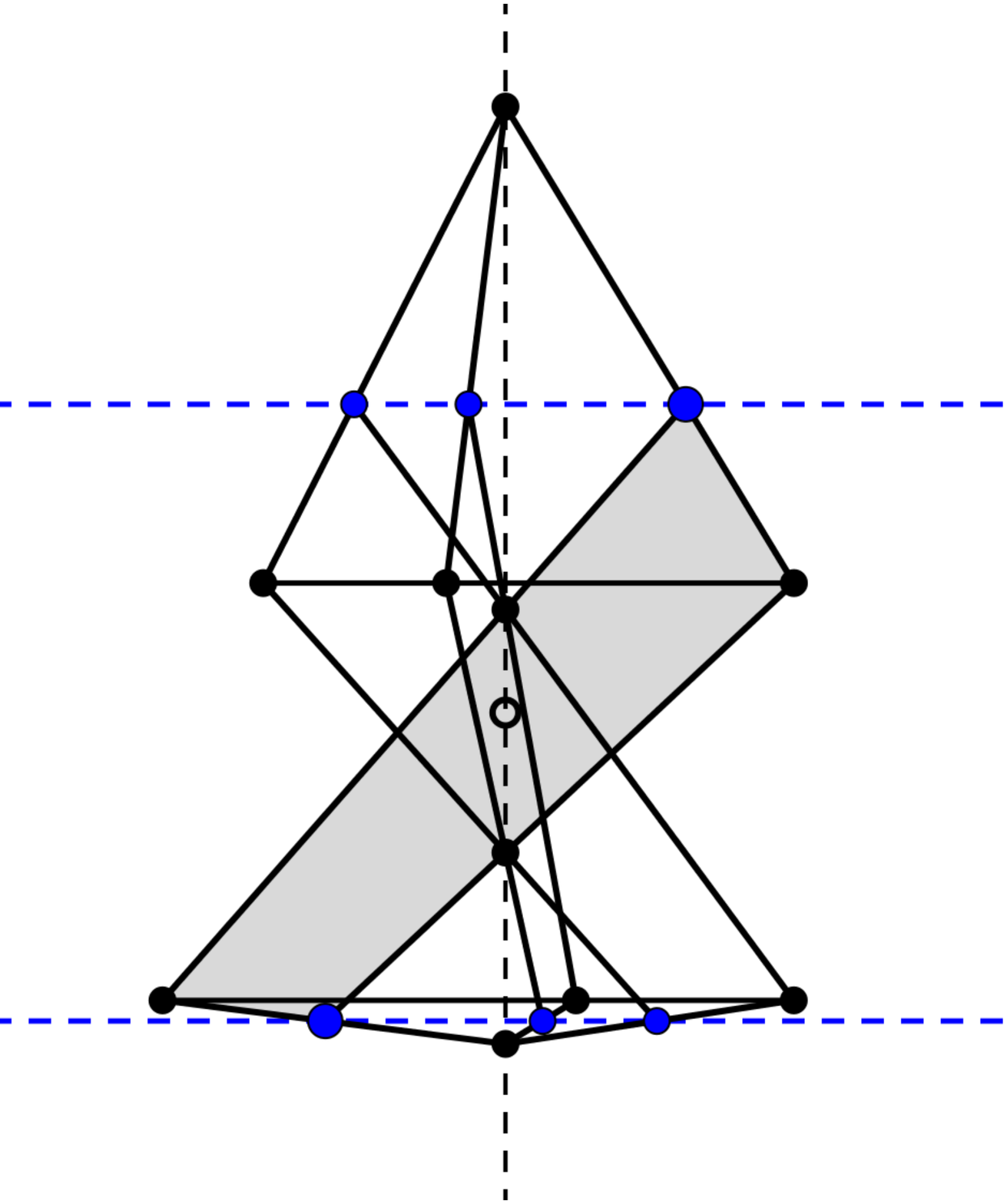}
    }
    \ffigbox{\caption{The complete configuration.}\label{fig:DU2end}}{
    \includegraphics[width=.3\textwidth]{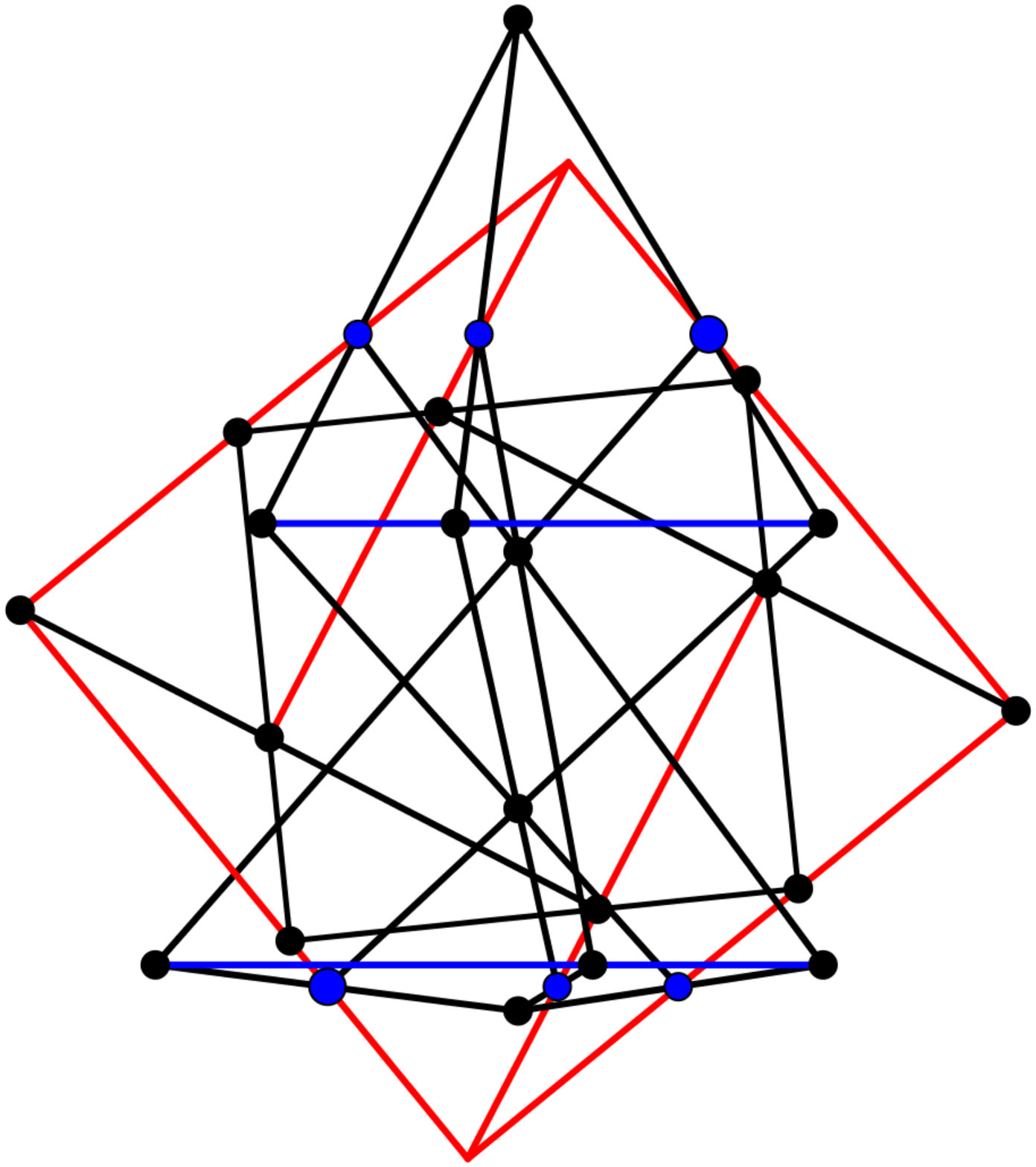}
      }
 \end{subfloatrow}   
 }   
{\caption{Example of the $\du^{(2)}(\mc C_{1}, \AR(\mc C_{2}))$ construction producing a $(26_3)$ configuration (for the details see the the text).}
\label{fig:DU_AR_Example}}
\end{center}
\end{figure}

Figure~\ref{fig:DU_AR_Example} shows an example of a $(26_4)$ configuration obtained by this construction. Here $\mc C_{1}$ is a $(12_3)$ configuration shown in Figure \ref{fig:DU2start}; 
its two pencils forming a centrally symmetric pair are highlighted in red, and their centre of symmetry is denoted by an empty circle. The points denoted by red empty 
circles will be deleted. In the $\AR(\mc C_{2})$ configuration shown in Figure \ref{fig:DU2middle}, $\mc C_{2}$ is a $(4_2)$ quadrilateral. The flexible lines are horizontal (the dashed copies will 
be deleted). The vertical dashed line is the axis of affinity. The starting quadrilateral is shaded; its two vertices forming a centrally symmetric pair are shown in a 
larger size. The completed configuration is shown in Figure \ref{fig:DU2end}.

\section{Systematic geometric constructions for $k = 5$}\label{sect:n5systematic}

\begin{figure}[htbp]
\begin{center}
\includegraphics[width = .48\linewidth]{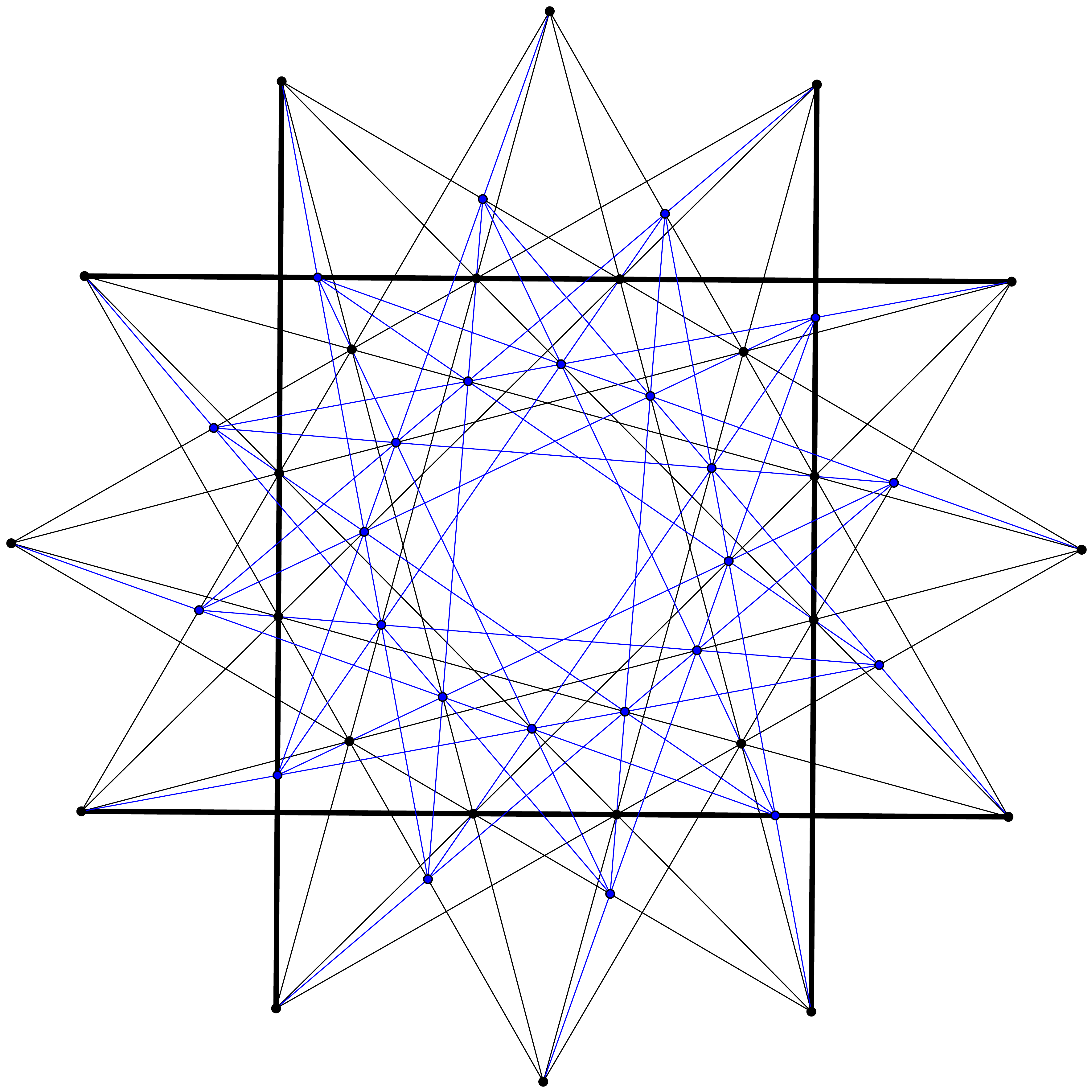} \hfill \includegraphics[width = .48\linewidth]{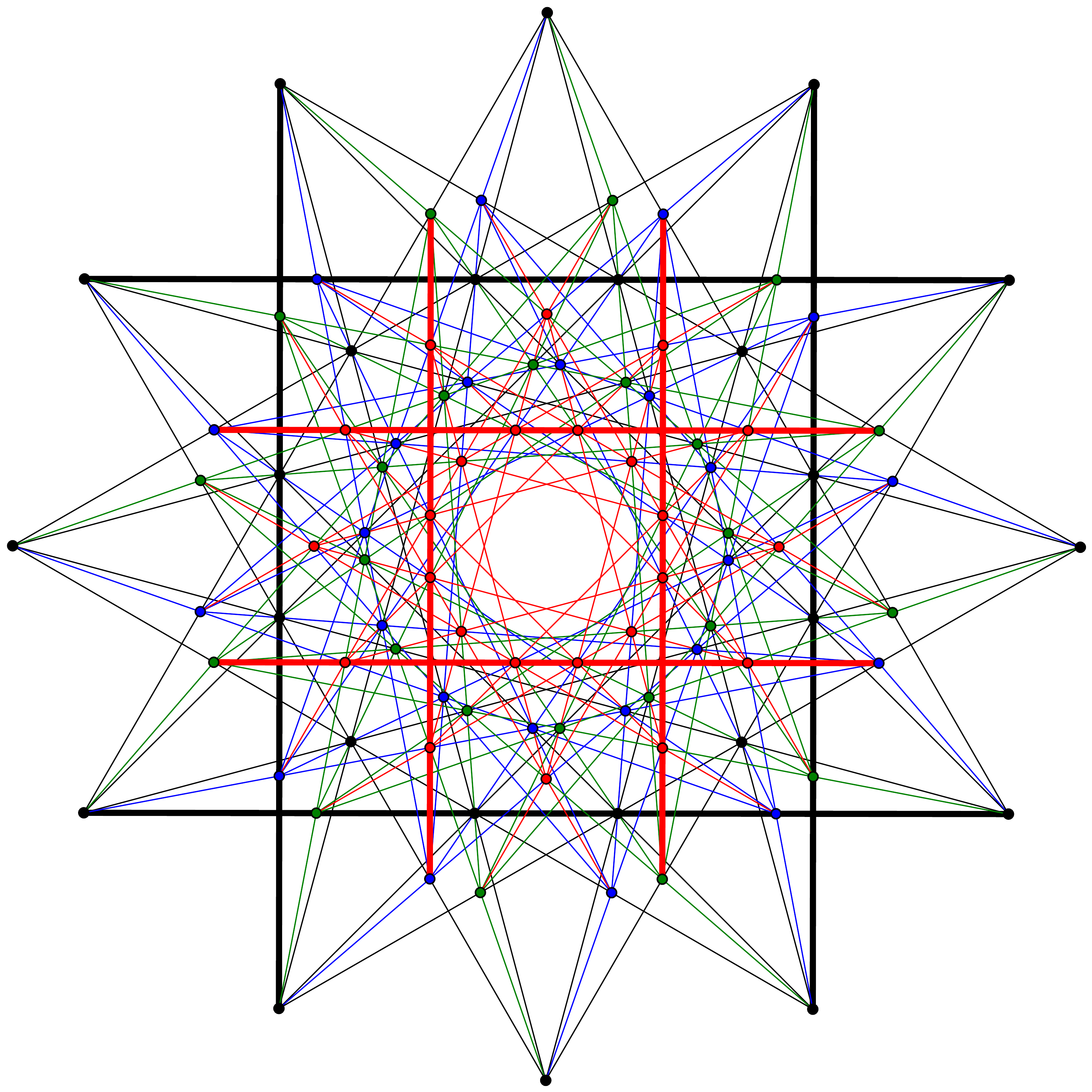}
\caption{The smallest known 5-configuration (left), a $(48_{5})$ configuration, and the smallest known 6-configuration (right), a $(96_{6})$ configuration. 
             Observe that the $(48_{5})$ configuration has two independent parallel pencils of size 2, and the $(96_{6})$ configuration has two independent parallel pencils of size 4.}
\label{fig:48-5-withParallels}
\end{center}
\end{figure}

The smallest known geometric 5-configuration is a $(48_{5})$ configuration (shown in Figure \ref{fig:48-5-withParallels}). It is unknown whether there are smaller 
geometric 5-configurations. In order to reduce the bound on $N_{5}$, in the next section we discuss various previously-known constructions that produce infinite series of 
(non-consecutive) 5-configurations, and introduce a few that have not been described previously in the literature. 

\subsection{``$\mc{A}$-series'' configurations}
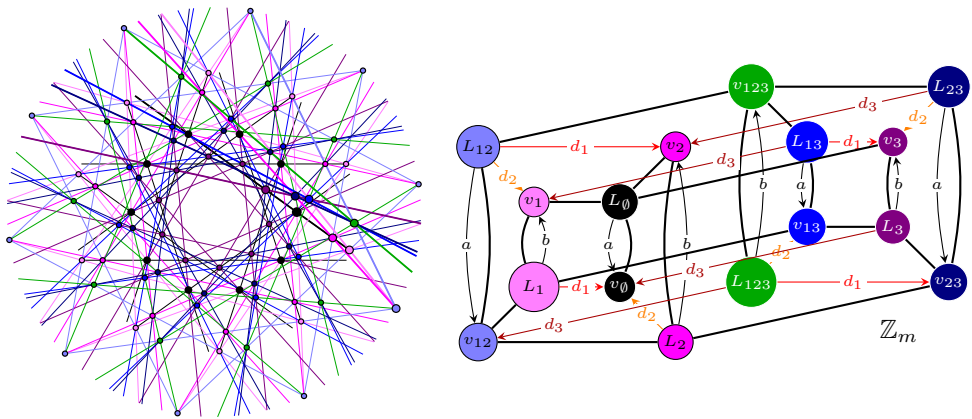
\begin{figure}[htbp]
\begin{center}
\ffigbox[]{
\begin{subfloatrow}[2]
\ffigbox{\caption{An example of a 5-configuration constructed as part of the $\mc{A}$-series, $\mc{A}(10;3,3;1,2,4)$}\label{fig:AseriesEx}}{
\includegraphics[width=\linewidth]{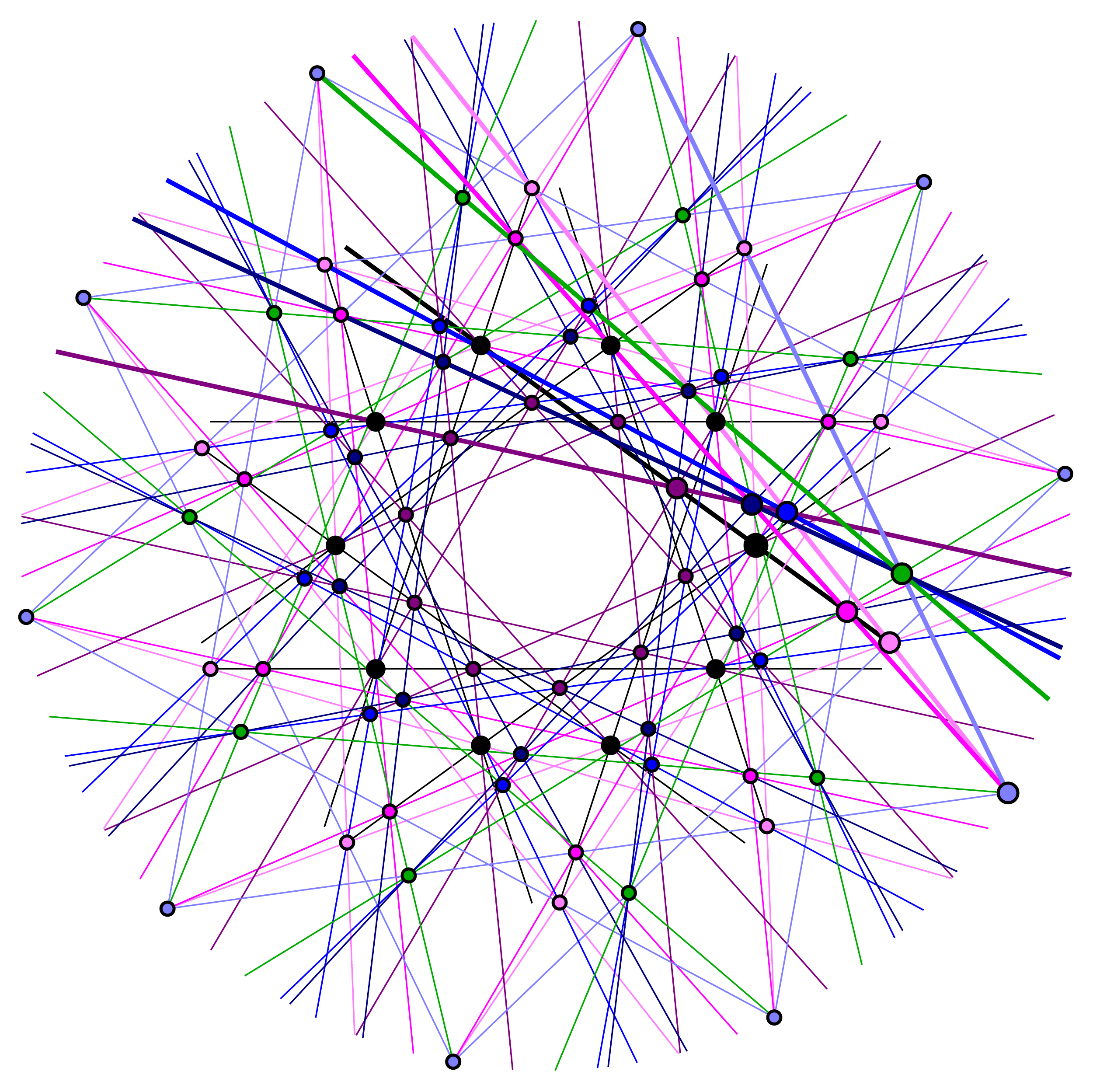}
}
\ffigbox{\caption{The reduced Levi graph for $\mc{A}(m;a,b;d_{1}, d_{2}, d_{3})$. Colors of symmetry classes match those in the configuration in (A). }\label{A-seriesRLG}}{
\begin{tikzpicture}[vtx/.style={draw, circle, font = \tiny, inner sep = 1 pt}, lin/.style={draw, circle, font = \tiny, inner sep = 1 pt}, lbl/.style={font = \tiny, inner sep = 1 pt, fill = white}, scale = .8
]
\definecolor{mag1}{rgb}{1,.5,1}
\definecolor{mag2}{rgb}{1,0,1}
\definecolor{mag3}{rgb}{.5,0,.5}
\definecolor{blu1}{rgb}{.5,.5,1}
\definecolor{blu2}{rgb}{0,0,1}
\definecolor{blu3}{rgb}{0,0,.5}
\definecolor{grn}{rgb}{0,.66,0}
\definecolor{drkRd}{rgb}{.66,0,0}

\node[vtx, white, fill=black,  ] (v-emp) at (-45:1){$v_{\emptyset}$};
\node[lin, white, fill=black, ] (L-emp) at (45:1){$L_{\emptyset}$};
\node[draw, circle, font = \tiny, inner sep = 1 pt, fill=mag1, ] (v1) at (45+90:1){$v_{1}$};
\node[draw, circle,font = \tiny, inner sep = 3 pt, fill=mag1, ] (L1) at (45+2*90:1){$L_{1}$};

\draw[<-] (v-emp) to[bend left=20] node[lbl]{$a$} (L-emp);
\draw[thick] (v-emp) to[bend right=20] (L-emp);

\draw[<-] (v1) to[bend left=20] node[lbl]{$b$} (L1);
\draw[thick] (v1) to[bend right=20] (L1);

\draw[thick] (L-emp) -- (v1);

\node[vtx, fill=mag2, ] (v2) at (45:2.3){$v_{2}$};
\node[lin, fill=mag2, ] (L2) at (-45:2.3){$L_{2}$};
\node[vtx, fill=blu1, ] (v12) at (-45-90:2.3){$v_{12}$};
\node[lin, fill=blu1, ] (L12) at (-45-2*90:2.3){$L_{12}$};

\draw[<-] (v2) to[bend left=10] node[lbl]{$b$} (L2);
\draw[thick] (v2) to[bend right=10] (L2);
\draw[<-] (v12) to[bend left=10] node[lbl]{$a$} (L12);
\draw[thick] (v12) to[bend right=10] (L12);

\draw[thick] (L-emp) -- (v2);
\draw[thick] (L1) -- (v12);
\draw[thick] (L2) -- (v12);

\tikzmath{\x1=4.5;\y1=1;}

\node[vtx, white, fill = mag3] (v3) at ($(L-emp)+(\x1,\y1)$) {$v_{3}$};
\node[vtx, white, fill = mag3] (L3) at ($(v-emp)+(\x1,\y1)$) {$L_{3}$};
\node[vtx, white, fill = blu2] (L13) at ($(v1)+(\x1,\y1)$) {$L_{13}$};
\node[vtx, white, fill = blu2] (v13) at ($(L1)+(\x1,\y1)$) {$v_{13}$};

\node[vtx, white, fill = blu3] (L23) at ($(v2)+(\x1,\y1)$) {$L_{23}$};
\node[vtx, white, fill = blu3] (v23) at ($(L2)+(\x1,\y1)$) {$v_{23}$};

\node[vtx, white, fill = grn] (L123) at ($(v12)+(\x1,\y1)$) {$L_{123}$};
\node[vtx, white, fill = grn] (v123) at ($(L12)+(\x1,\y1)$) {$v_{123}$};

\draw[<-] (v3) to[bend left=10] node[lbl]{$b$} (L3);
\draw[thick] (v3) to[bend right=10] (L3);
\draw[<-] (v13) to[bend left=10] node[lbl]{$a$} (L13);
\draw[thick] (v13) to[bend right=10] (L13);

\draw[<-] (v23) to[bend left=10] node[lbl]{$a$} (L23);
\draw[thick] (v23) to[bend right=10] (L23);
\draw[<-] (v123) to[bend left=10] node[lbl]{$b$} (L123);
\draw[thick] (v123) to[bend right=10] (L123);

\draw[->, red] (L1) -- node[lbl] {$d_{1}$} (v-emp);
\draw[->, red] (L12) -- node[lbl] {$d_{1}$} (v2);
\draw[->, red] (L13) -- node[lbl] {$d_{1}$} (v3);
\draw[->, red] (L123) -- node[lbl] {$d_{1}$} (v23);
\draw[->, orange] (L2) -- node[lbl] {$d_{2}$} (v-emp);
\draw[->, orange] (L12) -- node[lbl] {$d_{2}$} (v1);
\draw[->, orange] (L23) -- node[lbl] {$d_{2}$} (v3);
\draw[->, orange] (L123) -- node[lbl] {$d_{2}$} (v13);
\draw[->, drkRd] (L3) -- node[lbl, near end] {$d_{3}$} (v-emp);
\draw[->, drkRd] (L13) -- node[lbl, near start] {$d_{3}$} (v1);
\draw[->, drkRd] (L23) -- node[lbl, near start] {$d_{3}$} (v2);
\draw[->, drkRd] (L123) -- node[lbl, near end] {$d_{3}$} (v12);

\draw[thick] (L13) -- (v123);
\draw[thick] (L23) -- (v123);
\draw[thick] (L3) -- (v13);

\draw[thick] (L3) -- (v23);
\draw[thick] (L12) -- (v123);
\draw[thick] (L-emp) -- (v3);
\draw[thick] (L1) -- (v13);
\draw[thick] (L2) -- (v23);

\draw node[below left= .2 and .1 of v23] {$\mathbb{Z}_{m}$};

\end{tikzpicture}
}
\end{subfloatrow}

}{
\caption{The $(8m_{5})$ $\mc{A}$-series configurations, and their reduced Levi graphs}
\label{fig:Aseries}
}
\end{center}
\end{figure}

In \cite{BerFau2013}, one of the authors (LWB) and Jill Faudree described a construction method that produces $(n_{k})$ configurations for any arbitrary 
$k \geq 3$ (but where $n$ depends on the details of the construction and chosen parameters). For $k = 5$, the configurations $\mc{A}(m; 3,3; 1,2,4)$ 
(among many other 5-configurations) are guaranteed to exist for all $m \geq 7$ except $m = 12$ and produce configurations $(8m_{5})$. An example 
of the $(80_{5})$ configuration $\mc{A}(10;3,3;1,2,4)$ is shown in Figure \ref{fig:Aseries} along with the reduced Levi graph for the general family 
$\mc{A}(m; a,b;d_{1}, d_{2}, d_{3})$.  Moreover, a straightforward geometric argument shows that if $m$ is odd, we are guaranteed that for 
$\mc{A}(m; 3,3; 1,2,4)$,  $p = q = 2$, while if $m$ is even, $p = q = 4$. 

There are a few special cases among the $\mc{A}$-series configurations; of particular note is the configuration $\mc{A}(12; 4,4;1,3)$, 
which according to the details of the construction should be expected to produce a 4-configuration, but in fact has ``extra incidences'' 
and produces the smallest known 5-configuration, a $(48_{5})$ configuration. The configuration $\mc{A}(12;3,3;1,2,4)$ has other 
extra incidences that cause it to not be a 5-configuration. However, the $(96_{5})$ configuration $\mc{A}(12; 5,5;1,4,6)$ does exist and 
has $p = q = 4$. In Table \ref{table:5-cfg-constructions}, we use the notation $\mc{A}(m)$ to refer to the specific configurations $\mc{A}(m; 3,3;1,2,4)$ for $m \neq 12$.

\subsection{5-configurations that are $h$-astral in $\mathbb{E}^{+}$}

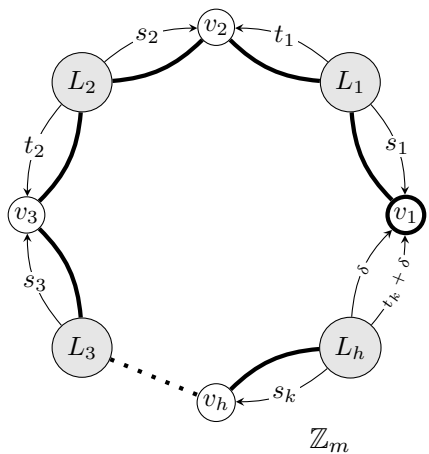
\begin{figure}[htbp]
\begin{center}

\begin{tikzpicture}[scale=0.5, ]
\tikzstyle{points}=[draw,circle,inner sep=1];
\tikzstyle{lines}=[draw,circle, fill=gray!20, inner sep=3];
\tikzstyle{points2}=[points,inner sep=1,font=\tiny];
\tikzstyle{lines2}=[lines,inner sep=4,font=\tiny];
\tikzstyle{arcleft}=[bend left=10,ultra thick];
\tikzstyle{arcright}=[bend right=10];
\tikzstyle{arclabel}=[fill=white,inner sep=1pt];

\def\x{45}
\def\rad{5}
\def\centerarc[#1](#2)(#3:#4:#5){ \draw[#1] ($(#2)+({#5*cos(#3)},{#5*sin(#3)})$) arc (#3:#4:#5); }

\node[points, ultra thick] (v00) at (0*\x:\rad) {$v_1$};
\node[lines ] (L00) at (1*\x:\rad) {$L_1$};
\node[points] (v01) at (2*\x:\rad) {$v_2$};
\node[lines ] (L01) at (3*\x:\rad) {$L_2$};
\node[points] (v02) at (4*\x:\rad) {$v_3$};
\node[lines ] (L02) at (+5*\x:\rad) {$L_3$};
\node[points, fill = white] (v0n1) at (6*\x:\rad) {$v_{h}$};
\node[lines ] (L0n1) at (7*\x:\rad) {$L_{h}$};

\tikzstyle{arcleft}=[bend left=20,ultra thick];
\tikzstyle{arcright}=[bend right=20];
\tikzstyle{arclabel}=[fill=white,inner sep=1];

\draw[arcright, <-] (v00) to node [arclabel] {$s_1$} (L00); \draw [arcleft] (v00) to (L00);
\draw[arcright, ->] (L00) to node [arclabel] {$t_1$} (v01); \draw [arcleft] (L00) to (v01);
\draw[arcright, <-] (v01) to node [arclabel] {$s_2$} (L01); \draw [arcleft] (v01) to (L01);
\draw[arcright, ->] (L01) to node [arclabel] {$t_2$} (v02); \draw [arcleft] (L01) to (v02);
\draw[arcright, <-] (v02) to node [arclabel] {$s_3$} (L02); \draw [arcleft] (v02) to (L02);

\draw[arcright, <-] (v0n1) to node [arclabel] {$s_{k}$} (L0n1); \draw [arcleft] (v0n1) to (L0n1);
\draw[arcright, ->] (L0n1) to node [arclabel, rotate=45+22.5
] {\tiny${t_k + \delta}$} (v00); \draw[arcleft, thin, ->] (L0n1) to node [arclabel, rotate=72] {\tiny$\delta$}(v00);


\begin{scope}[on background layer]
\draw[loosely dotted, ultra thick] (L02) -- (v0n1);

\end{scope}

\draw ($(v0n1) + (3,-1)$) node{\large$\mathbb{Z}_{m}$};
\end{tikzpicture}

\caption{The reduced Levi graph for a general $h$-celestial 4-configuration with symbol $m\#(s_{1}, t_{1}; s_{2}, t_{2}; \ldots; s_{h}, t_{h})$. The quantity $\delta = \frac{1}{2} \sum_{i=1}^{h}(s_{i} -t_{i})$ is the ``twist'' of the configuration.}
\label{fig:celestialRLG}
\end{center}
\end{figure}

A \emph{celestial} 4-configuration with symbol $m\#(s_{1}, t_{1}; \ldots; s_{h}, t_{h})$ is a polycyclic (actually, polydihedral) $(mh_{4})$ configuration with $h$ symmetry classes of points and lines under the action of $\mathbb{Z}_{m}$ and the reduced Levi graph shown in Figure \ref{fig:celestialRLG}. It has been described in a number of references; see, e.g., \cite[Sections 3.5--3.8]{Gru2009b} (under the name $h$-astral, although $h$-astral typically refers to any polycyclic configuration with $h$ symmetry classes, not just celestial ones) and \cite{BerBer2014}. There are several constraints that need to be satisfied for the parameters $s_{i}, t_{i}$ to correspond to a realizable geometric configuration, described the cited references, but for most of the examples in this paper, the parameters are \emph{trivial}, which means that the constraints are automatically satisfied.

Following \cite{BerGru2010}, we say that a \emph{diameter} of a celestial configuration passes through the center of the configuration and one point of the symmetry class $v_{0}$, and a \emph{mid-diameter} is the rotation of a diameter by $\frac{\pi}{m}$. If $m$ is even, then diameters connect points $(v_{0})_{i}$ and $(v_{0})_{i + \frac{m}{2}}$, and mid-diameters may or may not pass through configuration points depending on parity considerations of the celestial configuration symbol, while if $m$ is odd, mid-diameters and diameters coincide.

In \cite[p. 235]{Gru2009b}, Gr\"{u}nbaum describes two configurations which are $h$-astral in the extended Euclidean plane, 
formed by adding diameters to certain celestial 4-configurations to form pencils of 5 parallel lines, and then adding certain points 
at infinity at the intersections of the pencils in such a way that each point of the configuration lies on 5 lines, and each line passes 
through 5 points. He provides two examples of this construction: one produces a $(60_{5})$ configuration (which is actually 
\emph{astral}---that is, it has $\lfloor \frac{k+1}{2} \rfloor$ symmetry classes of points and lines, where 
$k = 5$---in $\mathbb{E}^{+}$) by adding diameters and points at infinity to $(2) 12\#(4,1;4,5)$, 
and one which produces a $(50_{5})$ configuration with 5 symmetry classes in $\mathbb{E}^{+}$, 
by adding diameters and points at infinity to the trivial celestial configuration $10\#(1,2;3,4;2,1;4,3)$.

\begin{figure}[htbp]
\begin{center}
\ffigbox[]{
\begin{subfloatrow}[2]
\ffigbox{\caption{The $(50_{5})$ configuration $\mc{D}_{4}(5)$, which can be notated $2\cdot 5\#(2,1;4,3;1,2;3,4)D\infty;MD\infty$.}\label{fig:celDiams1}}{
\includegraphics[width=\linewidth]{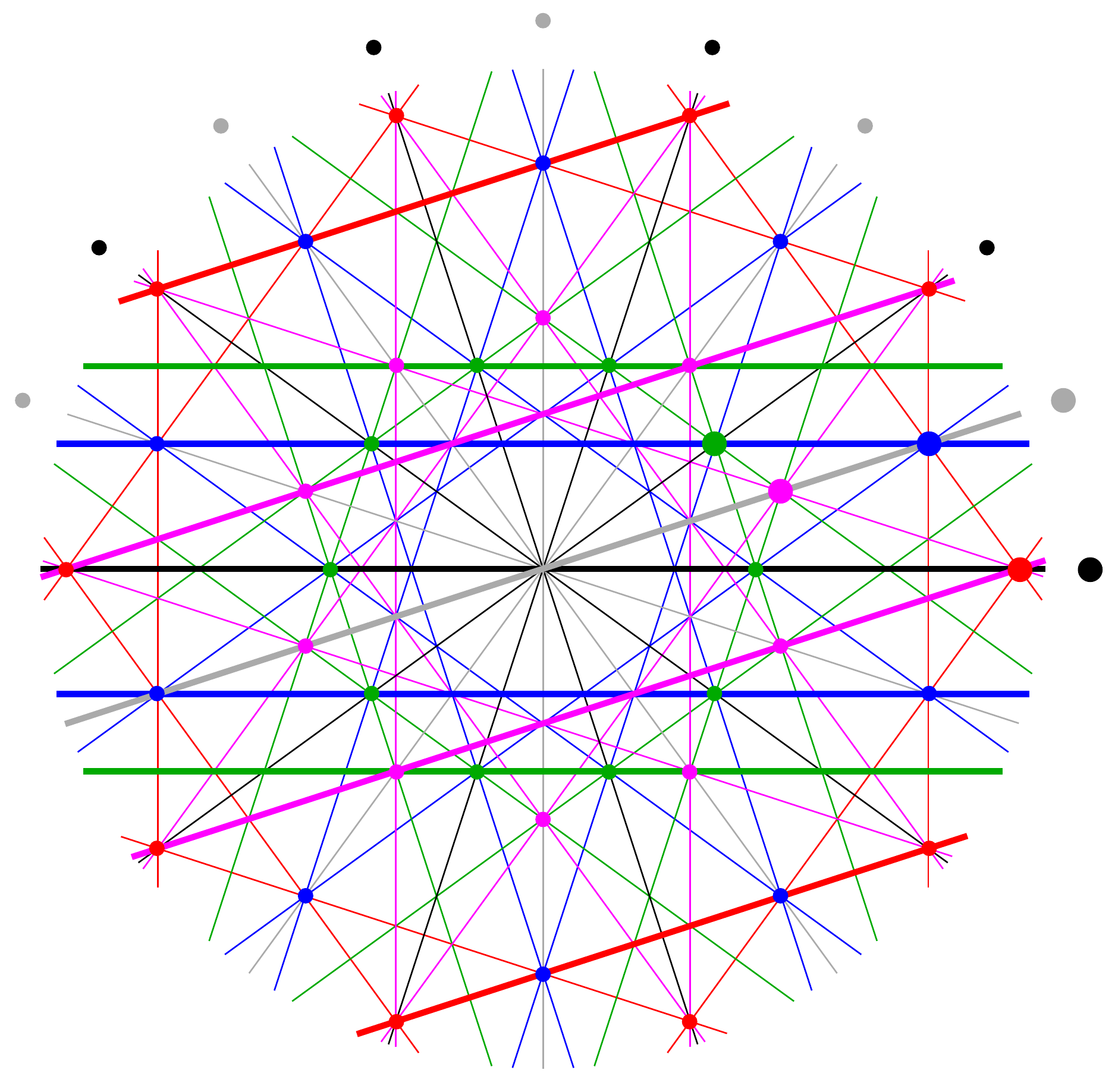}
}
\ffigbox{

\caption{A $(66_{5})$ configuration, $\mc{D}_{5}(5)$. One parallel pencil (in $\mathbb{E}^{2+}$) is shown with thick dashed lines.}
\label{fig:d5-11}
}{

\includegraphics[width = \linewidth]{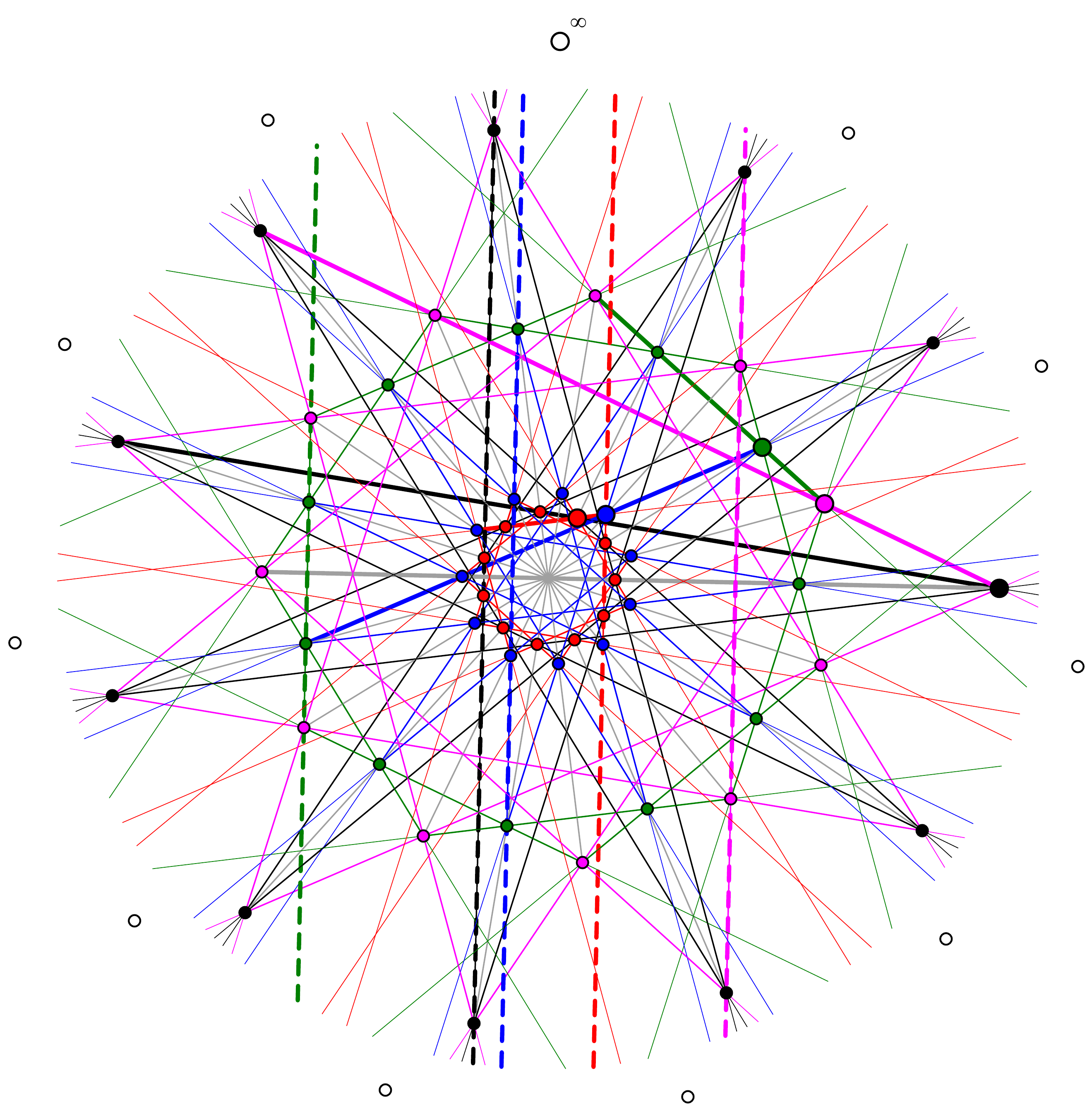}

}

\end{subfloatrow}
}{
\caption{
Examples of the $\mc{D}_{4}(\ell)$ and $\mc{D}_{5}(\ell)$ construction in $\mathbb{E}^{2+}$. Points outside the circular boundary of the configuration are at infinity.
}
\label{fig:D4-D5}
}
\end{center}
\end{figure}

Our first construction generalizes Gr\"{u}nbaum's $(50_{5})$ construction to form  $(10\ell_{5})$ configurations for $\ell \geq 1$, and our new second construction provides a  series of $(6(2\ell+1)_{5})$ configurations for $\ell \geq 5$. Small examples of each construction are shown in Figure \ref{fig:D4-D5}.

\subsubsection{The $(10\ell_{5})$ configurations $\mc{D}_{4}(\ell)$ based on the celestial $2\ell\#(2,1;4,3;1,2;3,4)$ configuration.}

\begin{construction}The $\mc{D}_{4}(\ell)$ construction proceeds as follows:
\be
\item Construct the trivial celestial configuration $2\ell\#(2,1,4,3,1,2,3,4)$, which has 4 symmetry classes of points and lines.
\item Construct all the diameters and mid-diameters (see \cite{BerGru2010} for discussion).
\item Add points at infinity in the directions of each of the diameters and mid-diameters.
\item Use a suitable projective transformation to project the configuration into $\mathbb{E}^{2}$.
\ee
\end{construction}

\begin{lem} The 
$\mc{D}_{4}(\ell)$
construction produces a series of $(10\ell_{5})$ configurations, for $\ell \geq 5$, with $p = q = 1$ when the configuration is projected into $\mathbb{E}^{2}$ from $\mathbb{E}^{2+}$.
\end{lem}

\begin{proof}
First, note that $2\ell\#(2,1;4,3;1,2;3,4)$ is only defined for $\ell \geq 5$. By considering the parity of the spans of the configuration 
(following the similar analysis from \cite{BerGru2010}), note that the red and green points lie on lines through the origin at each 
even multiple of $\frac{ \pi}{2\ell}$ (``diameters'' in the language of  \cite{BerGru2010}), while the blue and magenta points lie on lines 
at an odd multiple of $\frac{ \pi}{2\ell}$ (``mid-diameters'' in the language of  \cite{BerGru2010}). That is, each diameter and 
mid-diameter passes through four points, two each of two colors.

When $\ell$ is odd, the configuration has $\ell$ quadruples of blue-green parallels, parallel to the diameters, and
and $\ell$ quadruples of red-magenta parallels, parallel to each of the mid-diameters, for $2\ell$ sets of parallels, and by this description, each is also parallel to one of the diameters, 
that is, lines through the origin at an angle of $\frac{i\pi}{2\ell}$ for each $i$ with $0 \leq i < 2 \ell$. When $\ell$ is even, the blue-green parallels are parallel to the mid-diameters and the red-magenta parallels are parallel to the diameters.

In total, there are $(4\cdot 2\ell)$ points and $(4\cdot 2\ell)$ lines from the original configuration; $2\ell$ 
added diameters, and $2\ell$ added points at infinity, for a total of $5\cdot 2\ell$ points, each lying on 
5 lines, and $5\cdot 2\ell$ lines, each passing through 5 points.
\end{proof}

\subsubsection{The $((12\ell+6)_{5})$ configurations $\mc{D}_{5}(\ell)$ based on the 5-celestial configuration \\$(2\ell+1)\#(5,1;2,3;4,5;1,2;3,4)$.}

A second, new, construction, called $\mc{D}_{5}(\ell)$, begins with a trivial 5-celestial 4-configuration with pencils of 5 parallel lines and again adds diameters and points at infinity. In this case, the underlying celestial configurations have an odd number of points in each symmetry class.

\begin{construction} The $\mc{D}_{5}(\ell)$ construction proceeds as follows:
\be
\item Construct the trivial celestial configuration $(2\ell+1)\#(5,1;2,3;4,5;1,2;3,4)$ for any $\ell \geq 5$.
\item Add all diameters (which, since $2\ell +1$ is odd, are also mid-diameters).
\item Add points at infinity in each of the directions $\frac{2 i \pi}{2\ell+1}$.
\item If you like, use a suitable projective transformation to project the configuration into $\mathbb{E}^{2}$.
\ee
\end{construction}

\begin{lem} The construction 
$\mc{D}_{5}(\ell)$
produces a $(6(2\ell+1)_{5}) = ((12\ell+6)_{5})$ configuration for all $\ell \geq 5$, with $p = q = 1$ if the configuration is projected into $\mathbb{E}^{2}$.
\end{lem}

\begin{proof} 
First, observe that the configuration $(2\ell+1)\#(5,1;2,3;4,5;1,2;3,4)$ is only defined for $\ell \geq 5$. 
Following the arguments in \cite{BerGru2010} (applied to the situation where $m = 2q+1$ is odd), observe that diameters pass through one point of each of 
the five symmetry classes of points. Specifically, the diameter $D_{i}$, 
which passes through $(v_{0})_{i}$ (black) and the center of the configuration, also passes through $(v_{1})_{i-2}$ (red), $(v_{2})_{i+4}$, $(v_{3})_{i-1}$, and $(v_{4})_{i+5}$.

Next observe that the lines $(L_{0})_{i+3}$ (black), $(L_{1})_{i-3}$ (red), $(L_{2})_{i+2}$ (blue), $(L_{3})_{i+4}$ (green), $(L_{4})_{i-2}$ (magenta) are all parallel; we add a point at infinity, $\infty_{i}$, in each of those parallel directions.

Projecting via an appropriate projective transformation yields a 5-configuration in the ordinary Euclidean plane. Any parallels in the original configuration are likely destroyed in the projected version, but we certainly can find two independent lines in two different directions, hence $p = q = 1$.
\end{proof}

The smallest member of this family, a $(66_{5})$ configuration named $\mc{D}_{5}(11)$, is shown in Figure \ref{fig:d5-11}.

\subsection{Nesting celestial 4-configurations to produce 5-configurations}

In~\cite{BerJacVer2016}, new 5-configurations (in particular, the first known class of movable 5-configurations) were developed by ``nesting'' certain celestial 4-configurations 
and connecting them via repeated applications of two geometric lemmas, the Crossing Spans Lemma (see, e.g.,~\cite{BerFauPis2016}) and the Configuration Construction Lemma. 
In particular, given any $h$-celestial cohort $m\#S;T$ where $S = \{s_{1}, \ldots, s_{h}\}$ and $T = \{t_{1}, \ldots, t_{h}\}$, if $S \cap T = \varnothing$, then the construction in \cite{BerJacVer2016} produces a $((mh^{2})_{5})$ configuration. In what follows, we abbreviate certain of these configurations (with fixed choices of parameters) as $\mc{N}(h, \nicefrac{m}{s}, s)$, where the particular choices of $S$ and $T$ depend on $s$ and are described below.

While a complete understanding of parameters that correspond to celestial configurations for large values of $h$ does not exist, there are known results for $h = 2, 3, 4$ (see \cite{BerBer2014,Ber2001, Gru2009b} among others) that are particularly relevant for the construction of small $(n_{5})$ configurations.

\begin{prop} The following cohorts of systematic celestial configurations produce infinite classes of 5-configurations through the ``nesting'' construction with useful smallest elements:
\be
\item The systematic 2-celestial cohort $6s\#\{3s-1, 1\};\{2s, 3s-2\}$ produces a $(24s_{5})$ configuration with symbol $\mc{N}(2,6,s) = [6s\#\{3s-1, 1\};\{2s, 3s-2\}]*2$ for all $s \geq 2$; the case where $s = 2$ produces the known $(48_{5})$ configuration, which is the smallest known geometric 5-configuration.\footnote{It is interesting to note that this particular configuration can be re-analyzed as coming from several different constructions; it is also a special case of the $\mc{A}$-series construction, with symbol $\mc{A}(12; 3,3;1,4)$ (expected to be a 4-configuration), due to extra incidences.} Note that all elements of the celestial configuration $6s\#(3s-1, 2s; 1, 3s-2)$ have two pairs of disjoint parallel pairs, inherited by the nested configuration, so $p = q = 2$.

\item The systematic 3-celestial cohort $2s\#\{s-1,1,s-4\};\{s-2,2,s-2\}$ produces a $(3^{2}(2s)_{5}) = (18s_{5})$ configuration with symbol $\mc{N}(3,2,s) = [2s\#\{s-1,1,s-4\};\{s-2,2,s-2\}]*3$ for $s \geq 5$, $3 \nmid s$. The smallest element is a $(90_{5})$ configuration.
When $s$ is odd, the celestial configuration $2s\#(s - 1, s - 2; 1, s - 2; s - 4, 2)$ has two disjoint bands of 6 parallel lines inherited by the nested configuration. When $s$ is even, each of the three celestial configurations involved in the ``nesting'' process has one pencil of four parallel lines and one pair of two parallel lines, so using one pencil of four parallel lines from two of the nested configurations maximizes the values of $p$ and $q$. In summary, if $s$ is odd, $p= q = 6$, and if $s$ is even, $p = q = 4$.

\item  The systematic 3-celestial cohort $3s\#\{s+1, s-1, 1\};\{s,s,3\}$ produces a $(3^{2}(3s)_{5}) = (27s_{5})$ configuration with symbol $\mc{N}(3,3,s) = [3s\#\{s+1, s-1, 1\};\{s,s,3\}]*3$ for all $s \geq 3$ except $s= 4$; the case where $s = 3$ produces the known $(81_{5})$ configuration. When $s$ is even, each celestial configuration contains  4 parallels in each of two directions. When $s$ is odd and larger than $5$, they have disjoint triples of parallels in two directions. However, when $s = 3, 5$ in a single configuration, due to the small number of points of the celestial configuration resulting in non-disjoint pencils, there is only a triple of parallels. However, these can taken in each of two of the nested configurations. In summary, when $s$ is even, $p = q = 4$, and when $s$ is odd, $p = q = 3$.
\ee
\end{prop}

\begin{proof} These results follow directly from Algorithm 2 in \cite{BerJacVer2016} and the lists of systematic celestial configurations in \cite{BerBer2014}, along with ad hoc analysis of parallels in 2-celestial and 3-celestial configurations. 

Note that in the case $\mc{N}(3,2,s)$, the restriction that $3\nmid s$ is due to the fact that when $s$ is a multiple of 3, one of the permutations in the nested configurations results in celestial configurations with extra incidences (6 lines passing through one symmetry class of points, rather than 4).

In the case $\mc{N}(3,3,s)$, when $s = 4$, the celestial configuration $12\#(5, 3; 3, 4; 1, 4)$, which would be one of the participants in the nesting, is degenerate, and so must be excluded. \end{proof}

We introduce here a more restrictive variant of this construction, which involves nesting only two copies of a $h$-celestial 
configuration (producing a non-movable configuration); this construction produces different configurations from those formed by applying 
Algorithm 2 from \cite{BerJacVer2016} with $h= 2$, and the allowable cohorts are significantly more restrictive. We present the construction in generality, but in Section \ref{sect:bound5cfgs} we only use the $(54_{5})$ configuration constructed as $\mc{N}'(9) = [9\#\{4,2,1\}; \{3,3,3\}]*2'$, which is the third-smallest known 5-configuration.

\begin{construction} \label{const:newNest}
A construction for $[m\# \{s_{1}, s_{2}, s_{3}\}; \{t, t, t\}]*2'$ and \\
$[m\# \{s_{1}, s_{2}, s_{3}, s_{4}\}; \{t_{1}, t_{1}, t_{2}, t_{2}\}]*2'$, respectively, where $S \cap T = \varnothing$. 
\be
\item Construct the celestial 
configuration $m\#(s_{1},t; s_{2},t; s_{3},t)$ (respectively, \\$m\#(s_{1},t_{1}; s_{2},t_{2}; s_{3},t_{1}; s_{4}, t_{2})$)  with vertex classes $(v^{1}_{i})$ and $(L^{1}_{i})$ and centre 
$\mc{O}$, $i = 1,2,3$; (respectively, $i = 1,2,3,4$).
\item Choose a parameter $d$ so that the circle  passing through $(v^{1}_{3})_{d}$, $\mc{O}$ and $(v^{1}_{3})_{d-s_{2}}$ (respectively, $(v^{1}_{4})_{d-s_{2}}$) intersects line $(L^{1}_{1})_{0}$; there are typically several possibilities.
\item Construct the circle and choose an intersection with line $(L^{1}_{1})_{0}$ to form $(v^{2}_{1})_{0}$;
\item Using $(v^{2}_{1})_{0}$ as the initial 
starting vertex, construct the celestial configuration \\$m\#(s_{2},t; s_{3},t; s_{1},t)$ (respectively, $m\#(s_{2},t_{1}; s_{3},t_{2}; s_{4},t_{1}; s_{1},t_{2})$).  
\ee
\end{construction}

\begin{figure}[htbp]
\begin{center}
\ffigbox[]{
\begin{subfloatrow}[2]
\ffigbox{\caption{The general reduced Levi graph for $[m\# \{s_{1}, s_{2}, s_{3}\}; \{t, t, t\}]*2'$.}\label{fig:twiceNestedRLG}}{
\begin{tikzpicture}[scale = .75, vtx/.style = {draw, circle, font =\tiny, inner sep = 1 pt},  line/.style = {draw, circle, font =\tiny, inner sep = 2 pt}, lbl/.style = {font =\tiny, inner sep = 1pt, fill = white}]
\def\rad{3.5}
\node[vtx, fill=red] (v1) at (0:\rad) {$v^{1}_{1}$};
\node[line, fill=red, opacity=.5] (L1) at (360/6:\rad) {$L^{1}_{1}$};
\node[vtx, fill=orange] (v2) at (360*2/6:\rad) {$v^{1}_{2}$};
\node[line, fill=orange, opacity=.5] (L2) at (360*3/6:\rad) {$L^{1}_{2}$};
\node[vtx, fill=yellow] (v3) at (360*4/6:\rad) {$v^{1}_{3}$};
\node[line, fill=yellow, opacity=.5] (L3) at (360*5/6:\rad) {$L^{1}_{3}$};

\def\radd{2}
\def\angle{360/6+2*180/6}
\node[vtx, fill=blue] (v12) at (0+\angle:\radd) {$v^{2}_{1}$};
\node[line, fill=blue, opacity=.5] (L12) at (360/6+\angle:\radd) {$L^{2}_{1}$};
\node[vtx, fill=blue!30!white] (v22) at (360*2/6+\angle:\radd) {$v^{2}_{2}$};
\node[line, fill=blue!30!white, opacity=.5] (L22) at (360*3/6+\angle:\radd) {$L^{2}_{2}$};
\node[vtx, fill=mycyan] (v32) at (360*4/6+\angle:\radd) {$v^{2}_{3}$};
\node[line, fill=mycyan, opacity=.5] (L32) at (360*5/6+\angle:\radd) {$L^{2}_{3}$};

\darcr{v1}{s_{1}}{L1}
\darcl{v2}{t}{L1}
\darcr{v2}{s_{2}}{L2}
\darcl{v3}{t}{L2}
\darcr{v3}{s_{3}}{L3}
\draw[] (v1) edge [<-, bend left=20] node[midway, fill=white, inner sep=1 pt]{\tiny{$t+\delta$}} (L3);
\draw[] (L3)edge [bend left=20] node[midway, fill=white, inner sep=1 pt]{$\delta$}  (v1);

\darcr{v12}{s_{2}}{L12}
\darcl{v22}{t}{L12}
\darcr{v22}{s_{3}}{L22}
\darcl{v32}{t}{L22}
\darcr{v32}{s_{1}}{L32}
\draw[] (v12) edge [<-, bend left=20] node[midway, fill=white, inner sep=1 pt]{\tiny{$t+\delta$}} (L32);
\draw[] (L32)edge [bend left=20] node[midway, fill=white, inner sep=1 pt]{$\delta$}  (v12);

\draw[thick](L1) -- (v12);
\draw[thick, dashed] (L2) -- (v22);
\draw[thick, dashed] (L3) -- (v32);

\draw[thick] (L12) to node[midway, fill=white, inner sep=1 pt, lbl]{$d$}  (v3);
\draw[thick, dashed] (L22) to node[midway, fill=white, inner sep=1 pt]{\tiny{$d+\delta$}}  (v1);
\draw[thick, dashed] (L32) to node[midway, fill=white, inner sep=1 pt]{\tiny{$d+\delta$}} (v2);

\end{tikzpicture}
}
\ffigbox{\caption{The $(54_{5})$ configuration $\mc{N}'(9) = [9\#\{4,2,1\}; \{3,3,3\}]*2'$.}\label{[9{4,2,1}{3,3,3}*2']}}{
\includegraphics[width=\linewidth]{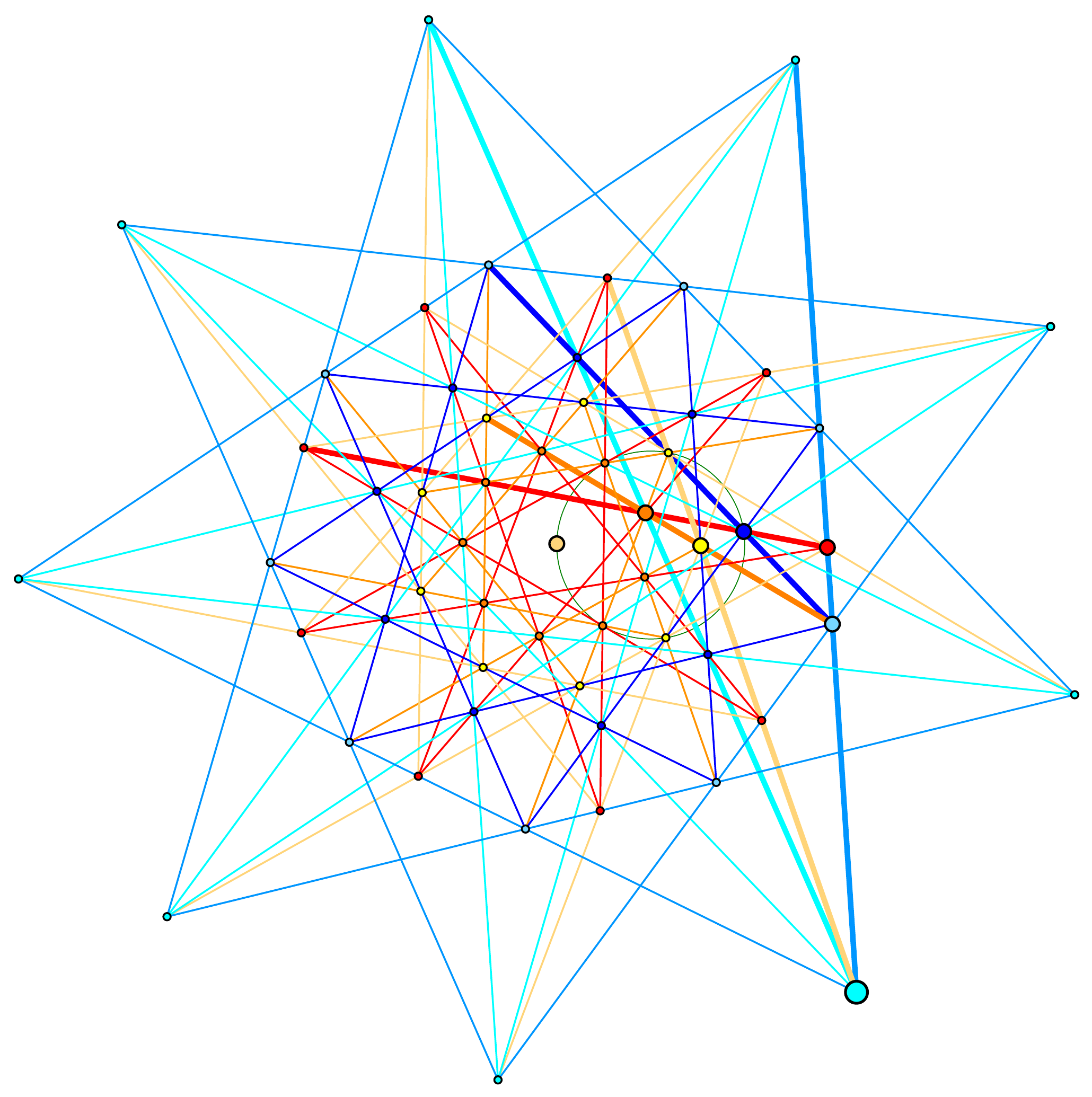}
}
\end{subfloatrow}

}{
\caption{A new class of $5$-configuration, denoted $[m\# \{s_{1}, s_{2}, s_{3}\}; \{t, t, t\}]*2'$, is constructed similarly 
to the construction in 
\cite{BerJacVer2016} 
and produces $(18m_{5})$ configurations in the rare cases when there exists 
a celestial 4-configuration of the cohort form $m\# \{s_{1}, s_{2}, s_{3}\}; \{t, t, t\}$. The smallest example, a $(54_{5})$ 
configuration, is $\mc{N}'(9) = [9\#\{4,2,1\}; \{3,3,3\}]*2'$, which has $p = 3$ and $q = 3$  (using a triple of parallels in each nested configuration). }
\label{fig:newNest}
}
\end{center}
\end{figure}
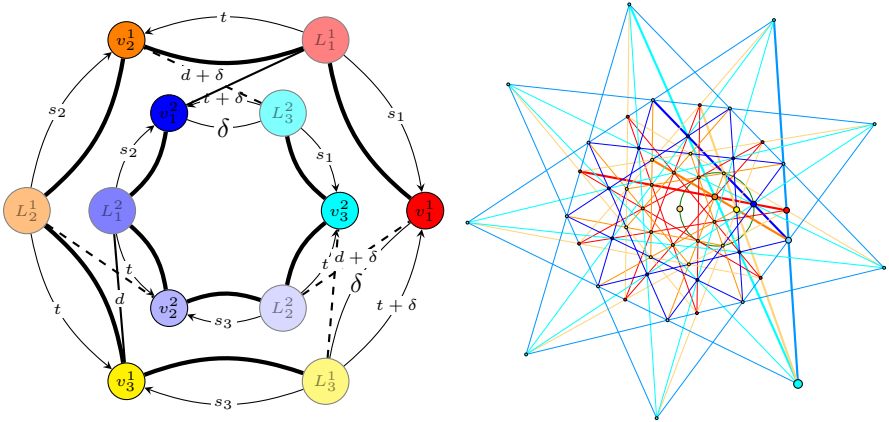

The reduced Levi graph of the construction for $[m\# \{s_{1}, s_{2}, s_{3}\}; \{t, t, t\}]*2'$ and the $(54_{5})$ configuration  $\mc{N}'(9) =[9\#\{4,2,1\}; \{3,3,3\}]*2'$ are shown in Figure \ref{fig:newNest}; as usual with celestial configurations, $\delta = \frac{1}{2}(\sum s_{i} - \sum t_{i})$. Analysis of the base celestial configuration shows that $p = q = 3$ in the $(54_{5})$ configuration (taking one pencil of 3 parallel lines in each sub-configuration, for example).

\begin{lem} Given a celestial cohort $m\# \{s_{1}, s_{2}, s_{3}\}; \{t, t, t\}$, Construction \ref{const:newNest} produces a $(6m_{5})$ configuration, and given a celestial cohort \\$m\# \{s_{1}, s_{2}, s_{3}, s_{4}\}; \{t_{1}, t_{1}, t_{2}, t_{2}\}$, Construction \ref{const:newNest} produces a $(8m_{5})$ configuration.
\end{lem}

\begin{proof} The argument to show that the resulting incidence structure is a 5-configuration is extremely similar to those given in the proof of 
Algorithm 2 from \cite{BerJacVer2016}, using the Crossing Spans Lemma and Configuration Construction Lemma described in that reference. Dotted lines in the reduced Levi graph shown in Figure \ref{fig:twiceNestedRLG} are forced by the Crossing Spans Lemma, and the highlighted green ``gadget'' in the reduced Levi graph is a result of the Configuration Construction Lemma, because in the construction, a circle was used to construct the points $(v^{2}_{1})_{0}$.
\end{proof}

This construction generalizes in a straightforward way.

\begin{lem}\label{lem:newNest} Given a celestial configuration cohort $m\# S; T$ with the property that $S \cap T = \varnothing$ and if $h$ is even, $T$ is partitioned into two subsets each of size $h/2$ where 
all the entries are the same, while if $h$ is odd, all elements of $T$ are equal, the construction $[m\#S;T]*2'$ produces a $(2mh_{5})$ configuration.
\end{lem}

Note that small celestial cohorts that satisfy the requirements of Lemma \ref{lem:newNest} are somewhat rare.  The construction $[15\#\{7,4,2,1\};\{6,6,3,3\}]*2'$ produces a $(120_{5})$ configuration using this new construction,  and $[24\#\{1,2,11\};\{8,8,8\}]*2'$ produces a $(144_{5})$ configuration, but we can produce these values of $n$ in other ways. There are no other cohorts with the necessary requirements for $m < 200$.%

It is also interesting to note that the infinite cohort 
$3q\#\{1, 2, \ldots, 2^{h-1}\};\{q, q, \ldots, q\}$ for $q = \frac{2^{h}+1}{3}$, $t$ 
odd, and $h >2$ described in \cite[Theorem 5.3]{BerJacVer2016} produces infinitely many cohorts of this type, including  
$33\#\{1, 2, 4, 8, 16\}, \{11, 11, 11, 11, 11\}$ which constructs a $(330_{5})$ configuration, but the size of 
the corresponding configurations grows quickly. This does let us conclude that there are infinitely many 6-astral 5-configurations in $\mathbb{E}^{2}$, however.

\section{Bound for 5-configurations}\label{sect:bound5cfgs}

In this section we establish a lower bound for the parameter $N_k$ defined in the Introduction, in the particular 
case of 5-configurations. 
The constructions used are summarized in Table \ref{table:5-cfg-constructions}.

\begin{table}[htp]
\caption{Summary of constructions for 5-configurations. By $\max_{1}(C,D)$ we mean the largest value among $\{p_{C}, p_{D}, q_{C}, q_{D}\}$ and $\max_{2}(C,D)$ is the second largest value among $\{p_{C}, p_{D}, q_{C}, q_{D}\}$, used with the \du\ construction, and $\max_{1}(C, r_{p}, r_{q})$ (resp.~$\max_{2}(C, r_{p}, r_{q})$) is the largest (resp.\ second largest) value from among $\{4(p_{C} - r_{p}), 4(p_{C} - r_{q}), 1\}$, used in the \AS\ construction.  Configurations marked with $*$ are line-flexible.}
\begin{center}

\begin{footnotesize}
\begin{tabular}{l p{.15\textwidth} l l  p{.2\textwidth}}
Construction & $n$ &$p$ & $q$ & conditions\\
\hline
\multicolumn{5}{l}{\emph{Systematic constructions}}\\
\hline
$\mc{A}(2s)$ & $16s$ & 4& 4 
& $s \geq 4$, $s\neq 6$\\
$\mc{A}(2s+1)$  & $16s+8$ & 2& 2& $s\geq 3$\\
$\mc{D}_{4}(s)$ & $10s$ & 1& 1 & $s \geq 5$ \\
$\mc{D}_{5}(s)$ & $12s+6$ & 1& 1 & $s \geq 5$ \\
$\mc{N}(2,6,s)$ 
& $24s$ & 2& 2& $s \geq 2$\\
$\mc{N}(3,2,6j+1)$
&$ 180j+18$ & 4& 4 & $j \geq 1$ \\ 
$\mc{N}(3,2,6j+2)$
&$ 180j+36$ & 6& 6 & $j \geq 1$ \\
$\mc{N}(3,2,6j+4)$
&$ 180j+72$ & 6& 6 & $j \geq 1$ \\
$\mc{N}(3,2,6j+5)$
&$ 180j+90$ & 4& 4 & $j \geq 0$ \\
$\mc{N}(3,3,2j+1)$
& $54j+27$ & 3 & 3 & $j \geq 1$\\ 
$\mc{N}(3,3,2j)$
& $54j$ & 4 & 4 & $j >2$\\ 
\hline
\multicolumn{5}{l}{\emph{Sporadic constructions}}\\
\hline
$\mc{N}'(9)$ & $54$ & 3&  3 
 \\ 
$\mc{A}(12)$ & 96 & 4 & 4 & $\mc{A}(12; 5,5; 1,2,4)$\\ 
\multicolumn{5}{l}{\emph{Gr\"unbaum Calculus constructions (* = flexible)}}\\
\hline
$\du(t)(C)$ & $(t+1)m-t$ &$p_{C}$&$q_{C}$ & $C = (m_{5})$\\
$\du(C,D^{*})$ & $m+\ell-1$ & $\max_{1}(C,D)$
 & $\max_{2}(C,D)$
& $C = (m_{5})$, $D = (\ell^{*}_{5})$\\
$\du^{(2)}(C, \AR(D))$ & $n+6m-2$ & $m-2$ & $\max\{p_{C}, q_{C}\}$& $C = (n_{5})$ with central symmetry; $D = (m_{4})$\\
$\AR(C)^{*}$ & $6m$ & $m$&$0$ & $C = (m_4)$ \\
$\PS(C)^{*}$ & $5m $ & $5 p_{C}$&$5 q_{C}$ & $C = (m_5)$ \\
$\AS(C, r_{p}, r_{q})^{*}$ & $4m+(p_{C} - r_p)+ (q_{C} - r_q)$ &$\max_{1}(C, r_{p}, r_{q})$& $\max_{2}(C, r_{p}, r_{q})$&  $C = (m_5)$, $1\leq r_{p} \leq p_{C}$, $0\leq r_{q} \leq q_{C}$\\
\end{tabular}
\end{footnotesize}
\end{center}
\label{table:5-cfg-constructions}
\end{table}%

\begin{thm}\label{thm:n5bound}
 $N_{5} \leq 166$. 
\end{thm}
\begin{proof}
In \cite{BerGevPis2021}, we showed that $N_{5}\leq 576$.

Furthermore, we have constructed a sequence of existing $(n_5)$ configurations from $m=166$ to $m=576$, with no gaps, as described below. 
This shows that $N_{5}\leq 166$. 

In the appendix, in Table \ref{table:n5cfgs}, we present a list of one example of an $(n_{5})$ configuration for each $n$ with $48\leq n \leq 250$, along with notations for missing configurations, choosing to include a configuration with a maximal known number of parallel lines. A list of all known configurations $48 \leq n \leq 576$ is available at \url{https://github.com/UP-LaTeR/5-bounds-6-bounds}. 

To confirm the correctness of the bounds in Theorem \ref{thm:n5bound} we performed two separate computer computations. The first computation used \emph{Mathematica} \cite{mathematica} and kept track of all construction information with each configuration constructed. The second approach used \emph{Sage} \cite{sagemath} and used the idea of combining arithmetic sequences of configurations.

The known $(n_{5})$ configurations were produced in \emph{Mathematica} as follows. At each step, the number of points and lines and the size of the independent pencils $p$ and $q$ for each configuration were tracked.
\be
\item Construct the numbers $n$ corresponding to the known systematic $(n_{5})$ configurations for $n =  48$ to 576 (note some numbers $n$ are produced from several systematic constructions):
\bi
\item $8m$ for $m = 7, ..., 72$ using $\mc A(m)$.
\item $10m$ for $m = 5,..., 57$ using $\mc{D}_{4}(s)$
\item $12s+6$ for $s = 5, ..., 47$ using $\mc{D}_{5}(s)$
\item $24s$ using  $\mc{N}(2,6,s)$; 
 in particular, this includes the known $(48_{5})$ configuration $\mc{N}(2,6,2)$
\item $18(3s + 2)$ and $18(3q + 4)$ for 
$q = 1, 2, ..., 7$
using  $\mc{N}(3,2,s) $
\item $27s$ for $s = 3, \ldots,32$ using $\mc{N}(3,3,s)$
\item Include $n = 54$ using the $(54_{5})$ configuration $\mc{N}'(9) = [9\#\{4,2,1\}; \{3,3,3\}]*2'$
\ei
\item Apply Affine Replication to the sequence of 4-configurations $\{ (18_{4})$, $(20_{4})$, $(21_{4})$, $(22_{4})$, $(24_{4})$,$(25_{4})$,$(26_{4})$, $\dotsc$, $(96_{4}) \}$ (consecutive after $n = 24$).

\item Next, apply each of the following constructions, at each step appending the results  with $n \leq 576$ to the list of known $(n_{5})$ configurations before applying the next step:
\be
\item Parallel Switch
\item $\du(t)$ for $t = 1, ..., 11$ (this produces lots of $(n_{5})$ configurations with $n>576$, since $\du(11)(
\mc{N}(2,6,2)
)$ produces a $(565_{5})$ configuration, but we ignore them)
\item Affine Switch
\item Parallel Switch
\item $\du(t)$ for $t = 1, ..., 11$ (Note that the single known configuration $(555_{5})$ is formed as $\du(1)(\PS(\mc{A}(7)))$, so we do need to do this second round).

\ee
\ee
Continuing to iterate the sequence of constructions (\PS, \du, \AS) more times does not result in any new configurations with $n \leq 576$. This sequence of constructions produces consecutive configurations for $n \geq 470$.

Next, observe that the outputs of Parallel Switch,  Affine Switch, and Affine Replication are all line-flexible configurations, so we use these outputs as one of the configurations $\mc{D}$ in the $\du(\mc C, \mc D^{*})$ flexible construction. This allows us to construct additional $(n_{5})$ configurations which were not constructible using the previous techniques. Specifically, we identify pairs of configurations $\mc C = (n_{5})$ and $\mc D = ((n_{x}^{*})_{5})$, where $x \in \{$PS, AR, AS$\}$ means that $D$ is the output of an application of parallel switch, affine replication, or affine switch respectively, such that $n+n_{x} - 1$ is missing from the list of previously constructed configurations. Then we apply $\du(\mc C, \mc D^{*})$ to those pairs to construct additional configurations. (Note that because the configurations $\mc D$ are line-flexible, we can use the \du\ construction to ``attach'' them to $\mc C$ without disrupting the parallel lines contained in $\mc C$.)

Finally, we find three small 5-configurations by hand using the $\du^{(2)}(C, \AR(D))$  construction, where $C$ is an $(n_{5})$ configuration with a halfturn 
symmetry relating two pencils and the configuration $D$ is a $(m_{4})$ configuration, as discussed 
in Section \ref{sec:DU-2}. Note that the $\mc{A}$-series configurations have halfturn symmetry for even $m$, and the known $(48_{5})$ configuration has 12-fold 
rotational symmetry (and thus halfturn symmetry) as well, so beginning with one point and pencil of the configuration, the image of the pencil under halfturn rotation is itself a pencil of the configuration, as required.

\begin{itemize}
\item $(182_{5}) = \du^{(2)}(\mc{A}(8), \AR(20_{4})) $ (note $\mc{A}(8)=(64_{5})$)
\item $(172_{5}) = \du^{(2)}(\mc{N}(2,6,2), \AR(21_{4})) $ (note $\mc{N}(2,6,2) = (48_{5})$)
\item $(154_{5}) = \du^{(2)}(\mc{N}(2,6,2), \AR(18_{4})) $ 
\end{itemize}

The final list of missing 5-configuration numbers greater than $48$ is $\{$49, 51, 52, 53, 55, 57, 58, 59, 61, 62, 63, 65, 67, 68, 69, 71, 73, \
74, 75, 76, 77, 79, 82, 83, 84, 85, 86, 87, 89, 91, 92, 93, 94, 97, \
98, 101, 103, 105, 106, 109, 113, 115, 116, 117, 118, 121, 122, 123, \
124, 125, 129, 133, 134, 137, 141, 145, 146, 147, 149, 151, 153,  \
158, 164, 165$\}$.
\end{proof}

It is unknown whether there exist $(n_{5})$ configurations with $n$ in the list of missing 5-configurations. Note that the smallest possible combinatorial $5$-configuration 
is $(21_{5})$ and the smallest \emph{known} topological 5-configuration is $(36_{5})$ (unpublished, found by the first author); the existence of geometric $(n_{5})$ 
configurations for $21 \leq n < 48$ is unknown as well.


\section{Bound for 6-configurations}\label{sect:bound6cfgs}


We use similar techniques to reduce the lower bound for $N_{6}$. In \cite{BerGevPis2021}, we showed that $N_{6} < 7350$, by using the fact that $N_{4}\leq 24$, so that there is a consecutive sequence of configurations $[(24_{4}), (25_{4}), (26_{4}), \ldots]$, using affine replication  to produce the sequence $[(144_{5}), (150_{5}), (6\cdot 26_{5}), \ldots]$, using affine replication again to produce the sequence $[(1008_{6}), (1050_{5}), (7\cdot 6\cdot 26_{5}), \ldots]$ and then applying affine switch to produce consecutive bands of 6-configurations for all $n \geq 7350$.

The smallest known 6-configuration is a $(96_{6})$ configuration,  constructed using the $\mc{A}$-series construction with parameters that lead to extra incidences in just the right way. Its construction is described in \cite{Ber2014} and it is shown in Figure \ref{fig:48-5-withParallels}.

To decrease the bound and fill in gaps between $n = 96$ and $n = 7350$, as in the 5-configuration case, we consider known constructions for systematic 6-configurations, as well as the Gr\"unbaum calculus operations described above. The three smallest known 6-configurations have $n = 96, 110, 112$, although the $(112_{6})$ configuration has symmetry classes of points that are very close together, so it is not very intelligible.

Known systematic constructions for 6-configurations are described in \cite{Ber2014}; there are two known infinite families. The first is the $\mc A$-series described above for 5-configurations, but using one more parameter: the family $\mc{A}(m; 3,3; 1,2,4,5)$ for $m\geq 7$ produces $(16m_{6})$ configurations, excluding $m = 12$ due to extra incidences. However, the configuration $\mc{A}(12; 5,5;1,2,3,7)$ is well-defined. In Table \ref{table:6-cfg-constructions}, we use the notation $\mc{A}(m)$ to refer to the specific configurations $\mc{A}(m; 3,3;1,2,4, 5)$ for $m \neq 12$. As with the $\mc{A}$-series for 5-configurations, these configurations have $p = q = 4$ when $m$ is even and $p = q = 2$ when $m$ is odd. The configuration $\mc{A}(7; 3,3;1,2,4,5)$ produces the known $(112_{6})$ configuration.

\begin{figure}[htbp]
\begin{center}
\ffigbox{
\begin{subfloatrow}[2]
\ffigbox{\caption{A $(110_{6})$ multicelestial configuration, $11\#(3,2,1)(4,5)$ }\label{fig:110-6}}{
\includegraphics[width=\linewidth]{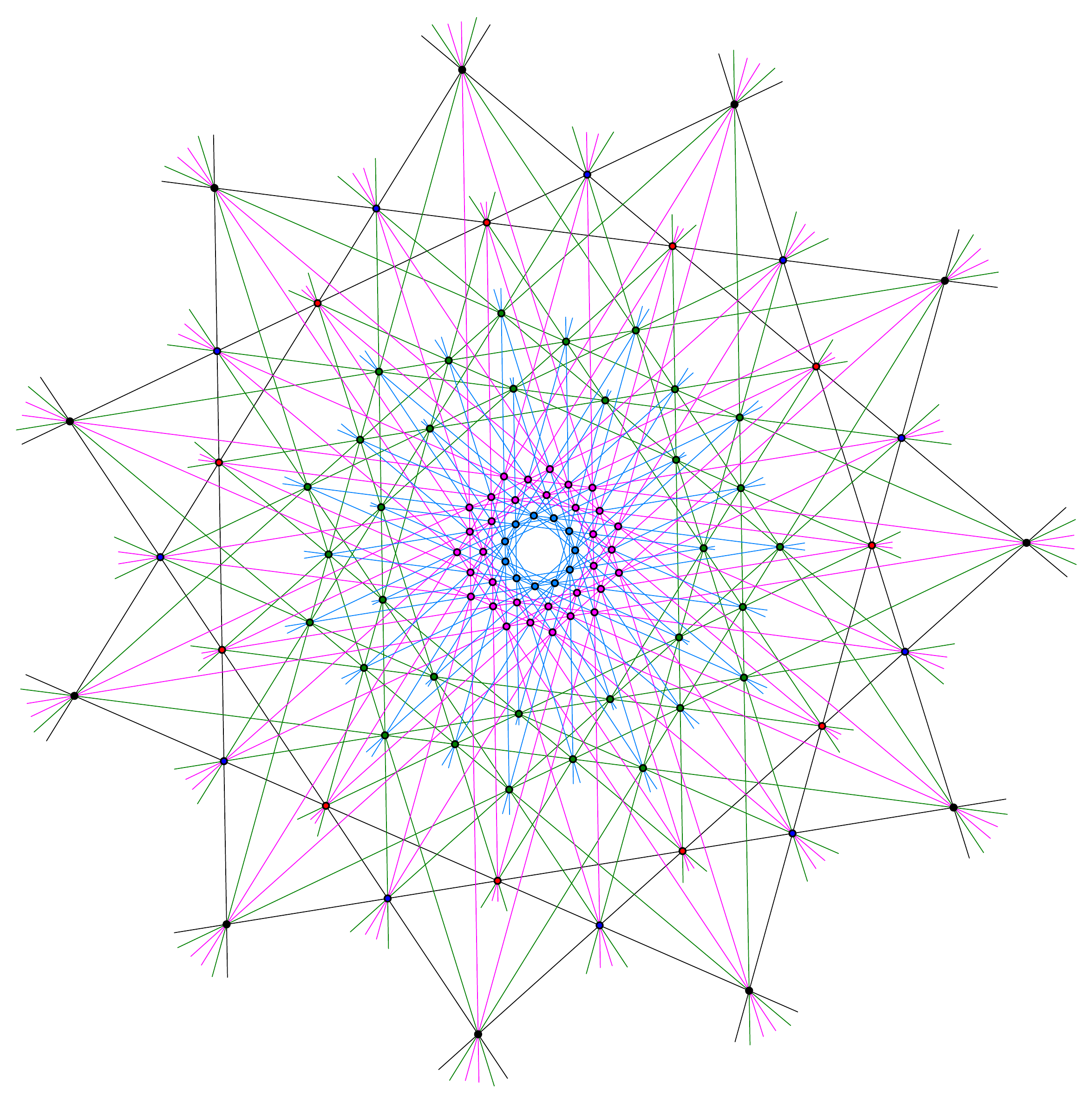}
}
\ffigbox{\caption{The $A$-series $(144_{6})$ configuration $\mc{A}(m;3,3;1,2,4,5)$ with $m = 9$. The configurations $(112_{6})$ and $(128_{6})$ configurations with $m = 7,8$ are well-defined but their points are too close together to be intelligible.}\label{fig:Aseries-6}}{
\includegraphics[width=\linewidth]{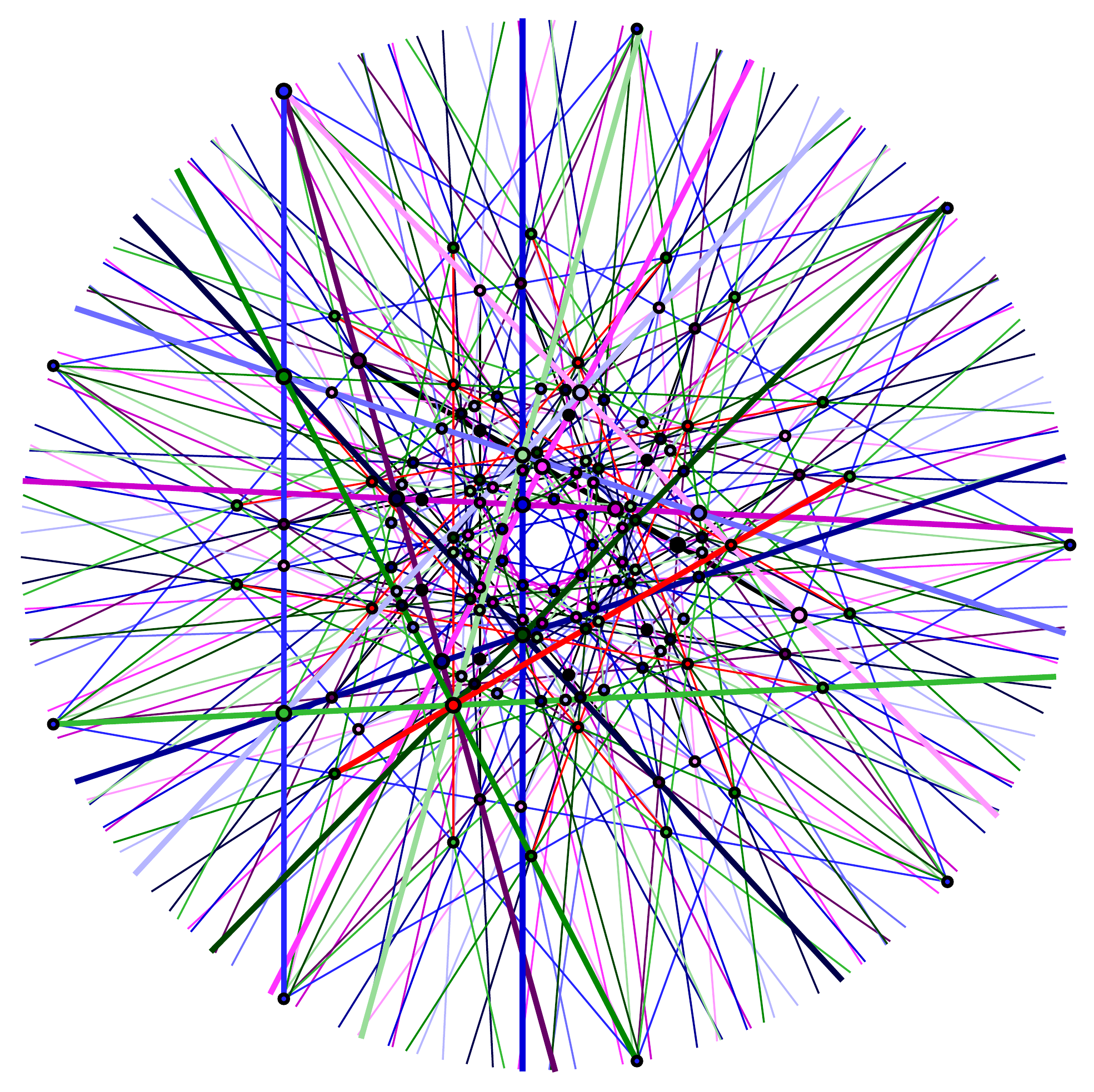}
}
\end{subfloatrow}
}{
\caption{Some pictures of small 6-configurations.}
\label{fig:6cfgs}
}
\end{center}
\end{figure}

The second known infinite family of 6-configurations, first described as part of a general family of $(2q,2k)$-configurations in \cite{Ber2011}, are called \emph{multicelestial configurations}. The construction of multicelestial 6-configurations with symbol  $m\#(t_{0}, t_{1}, t_{2})(s_{0},s_{1})$ is specifically detailed in \cite{Ber2014}. These produce $(10m_{6})$ configurations with dihedral $m$-gonal symmetry for all integers $m \geq 11$ except $m = 12$, with the known $(110_{6})$ configuration arising from $11\#(3,2,1)(4,5)$.  In general, these configurations have a lot of parallel lines, but the number of parallels depends on the choice of parameters. The parameters $m\#(3,2,1)(4,5)$ are valid for $m =  11$ and $m \geq 13$ ($m = 12$ has extra incidences), and the parameters $m\#(5,3,1)(7,9)$ are valid for $m \geq 19$. Note that while the multicelestial construction does not produce a $(120_{6})$ configuration, there are several known $(120_{6})$ configurations, using a variant construction, described in \cite{Ber2014} (e.g., the one shown in Figure 7 of that paper); inspection shows $p = q = 4$. Examples of a multicelestial configuration and an $\mc{A}$-series 6-configuration are shown in Figure \ref{fig:6cfgs}.

Because we will use multicelestial configurations as input into Gr\"unbaum Calculus operations, specifically into Affine Switch, it is useful to perform a detailed analysis of the maximum number of independent parallel lines in these configurations. Multicelestial configurations have 10 symmetry classes of points and 10 symmetry classes of lines, where the 10 line classes are indexed as $(L)_{i}$ (one class), $(L_{j}^{r})_{i}$ where $j \in \{0,1,2\}$ and $r \in \{0,1\}$ (6 classes), and $(L^{0,1}_{j,k})_{i}$ where $\{j,k\} \subseteq\{0,1,2\}$ (3 classes). The parameters $j$, $k$, $r$ refer to indices in the ordered sets $T = (t_{0}, t_{1}, t_{2})$ or $(s_{0}, s_{1})$ where $0 < t_{j}, s_{r} < \frac{m}{2}$ and $S \cap T = \varnothing$.  Analyzing the details of the construction from \cite{Ber2014} shows that lines $(L_{j}^{r})_{\frac{1}{2}(t_{j}-s_{r})+i}$  and lines $(L^{0,1}_{j,k})_{\frac{1}{2}(t_{j}+t_{k} - s_{0}-s_{1})+ i}$ are parallel to line $(L)_{0}$ whenever the indices (i.e., $\frac{1}{2}(t_{j}+t_{k} - s_{0}-s_{1})$ or $\frac{1}{2}(t_{j}-s_{r})$)  are integers. Furthermore, if $m$ is odd, we can replace $\frac{1}{2}$ with $\lceil \frac{m}{2} \rceil$ in all cases, so that for each of the other nine line types, there is one representative that is parallel to $(L)_{0}$. When $m$ is odd, some line types have representatives that are parallel to $(L)_{0}$ and some do not (depending whether the indices $\frac{1}{2}(t_{j}+t_{k} - s_{0}-s_{1})$ or $\frac{1}{2}(t_{j}-s_{r})$ are integers); however, it is true that if $(L^{y}_{x})_{a}$ is parallel to $(L)_{0}$, then so is  $(L^{y}_{x})_{a+\frac{m}{2}}$. 

In summary, if $m$ is odd, any value of parameters $(t_{0}, t_{1}, t_{2})(s_{0},s_{1})$ produces a pencil of 10 parallel lines in a single direction (say, parallel to $(L)_{0}$). However, when $m$ is even, analysis of the parities shows that if the parameters $(t_{0}, t_{1}, t_{2})(s_{0},s_{1})$ contains elements of both parities, then the best we can do is six lines in a single direction. Since $m\#(t_{0}, t_{1}, t_{2})(s_{0},s_{1})$ is disconnected when $m$ and the parameters are all even, it follows that the most parallel lines in a single direction for even $m$ occurs when $t_{j}$, $s_{r}$ are all odd, which is possible only when $m \geq 20$.

Finally, we need to determine whether the lines in the pencils are disjoint. For $m \geq 21$, $m$ odd, the configuration $m\#(5,3,1)(7,9)$ has two disjoint pencils with 10 lines in each pencil; similarly, for $m \geq 40$, $m$ even, the configuration $m\#(5,3,1)(7,9)$ has two disjoint pencils with 20 lines in each pencil. For smaller values of $m$, we inspected all configurations for each set of possible parameter values ($m\leq 20$) or odd parameter values ($22\leq m < 40$, $m$ even) after eliminating all lines that are incident with points of an initial pencil of lines parallel to $(L)_{0}$, to determine the maximum number of parallel lines in independent pencils. The results are shown in Table \ref{table:maxParallelsMulticelestial}.

\begin{table}[htp]
\caption{Maximum number of independent lines in two parallel pencils in a multicelestial 6-configuration $m\#(t_{0}, t_{1}, t_{2})(s_{0}, s_{1})$.}
\begin{center}
\begin{tabular}{l l l l}
odd $m$ & parameters & $p$ & $q$\\ \hline
11 & (3,2,1)(4,5) & 10 & 0\\
13& (3,2,1)(4,5) & 10 & 1\\
15& (3,2,1)(4,5) & 10 & 1 \\ %
17& (3,2,1)(4,5) & 10 & 2 \\
19& (5,3,1)(7,9) & 10& 4 \\
$\geq 21$& (5,3,1)(7,9) & 10& 10 \\ 
\end{tabular}
\begin{tabular}{l l l l}
even $m$ & parameters & $p$ & $q$\\ \hline
12 &- \\
14& (3,2,1)(4,5) & 12&0\\
16& (5,4,2)(6,7) & 12&4 \\ %
18& (5,4,1)(6,7) & 12&4 \\
20& (5,3,1)(7,9) & 20&0\\ 
22& (5,3,1)(7,9) & 20&0 \\
24& (5,3,1)(7,9) & 20&2 \\
26& (5,3,1)(7,9) & 20&2 \\
28& (5,3,1)(7,9) & 20&4 \\
30& (5,3,1)(7,9) & 20&2 \\
32& (5,3,1)(7,9) & 20&8 \\
34& (5,3,1)(7,9) & 20&4 \\
36& (5,3,1)(7,9) & 20&8 \\
38& (5,3,1)(7,9) & 20&8 \\
$\geq 40$& (5,3,1)(7,9) & 20&20 \\
\end{tabular}
\end{center}
\label{table:maxParallelsMulticelestial}
\end{table}%

\begin{table}[htp]
\caption{Summary of constructions for 6-configurations. By $\max_{1}(C,D)$ we mean the largest value among $\{p_{C}, p_{D}, q_{C}, q_{D}\}$ and $\max_{2}(C,D)$ is the second largest value among $\{p_{C}, p_{D}, q_{C}, q_{D}\}$, used with the \du\ construction, and $\max_{1}(C, r_{p}, r_{q})$ (resp. $\max_{2}(C, r_{p}, r_{q})$) is the largest (resp. second largest) value from among $\{5(p_{C} - r_{p}), 5(p_{C} - r_{q}), 1\}$, used in the \AS\ construction.  Configurations marked with $*$ are line-flexible.}
\begin{center}

\begin{footnotesize}
\begin{longtable}{l p{.15\textwidth} l l  p{.2\textwidth}}
Construction & $n$ &$p$ & $q$ & conditions\\
\hline
\multicolumn{5}{l}{\emph{Systematic constructions}}\\
\hline
$\mc{A}(2j+1)$ & $32j+16$ & 2 & 2 & $j \geq 3$\\
$\mc{A}(2j)$ & $32j$ & 4 & 4 & $j \geq 4$, $j\neq 6$\\
$m\#(t_{0}, t_{1}, t_{2})(s_{0}, s_{1})$ & $10m$ & see Table \ref{table:maxParallelsMulticelestial} & see Table \ref{table:maxParallelsMulticelestial} & $m \geq 11$, $m \neq 12$\\
\hline
\multicolumn{5}{l}{\emph{Sporadic constructions}}\\
\hline
$\mc{A}(12;4,4;1,3,5)$ & 96 & 4 & 4 & \\
$\mc{A}(12)$ & $192$ & 4 & 4 &  $\mc{A}(12; 5,5; 1,2,3,7)$\\
$4$-astral  & $120$ & 4& 4 
& See \cite[Figure 7]{Ber2014}\\
\hline
\multicolumn{5}{l}{\emph{Gr\"unbaum Calculus constructions (* = flexible)}}\\
\hline
$\du(t)(C)$ & $ (t+1)m-t$ &$p_{C}$&$q_{C}$ & $C = (m_{6})$\\
$\du(C,D^{*})$ & $m+\ell-1$ & $\max_{1}(C,D)$& $\max_{2}(C,D)$
& $C = (m_{6})$, $D^{*} = (\ell_{6})$\\
$\du^{(2)}(C, \AR(D))$ & see Table \ref{table:du2-6cfgs} &  $m-2$ & $\max\{p_{C}, q_{C}\}$ &  $C = (n_{6})$ with central symmetry; $D = (m_{5})$\\
$\AR(C)^{*}$ & $7m$ & $m$&$0$ & $C = (m_5)$ \\
$\PS(C)^{*}$ & $6m $ & $6 p_{C}$&$6 q_{C}$ & $C = (m_6)$ \\
$\AS(C, r_{p}, r_{q})^{*}$ & $5m+(p_{C} - r_p)+ (q_{C} - r_q)$ &$\max_{1}(C, r_{p}, r_{q})$& $\max_{2}(C, r_{p}, r_{q})$
 &  $C = (m_6)$, $1\leq r_{p} \leq p_{C}$, $0\leq r_{q} \leq q_{C}$\\
\end{longtable}
\end{footnotesize}
\end{center}
\label{table:6-cfg-constructions}
\end{table}%

\begin{figure}[htbp]
\begin{center}
\includegraphics[width = .7\linewidth]{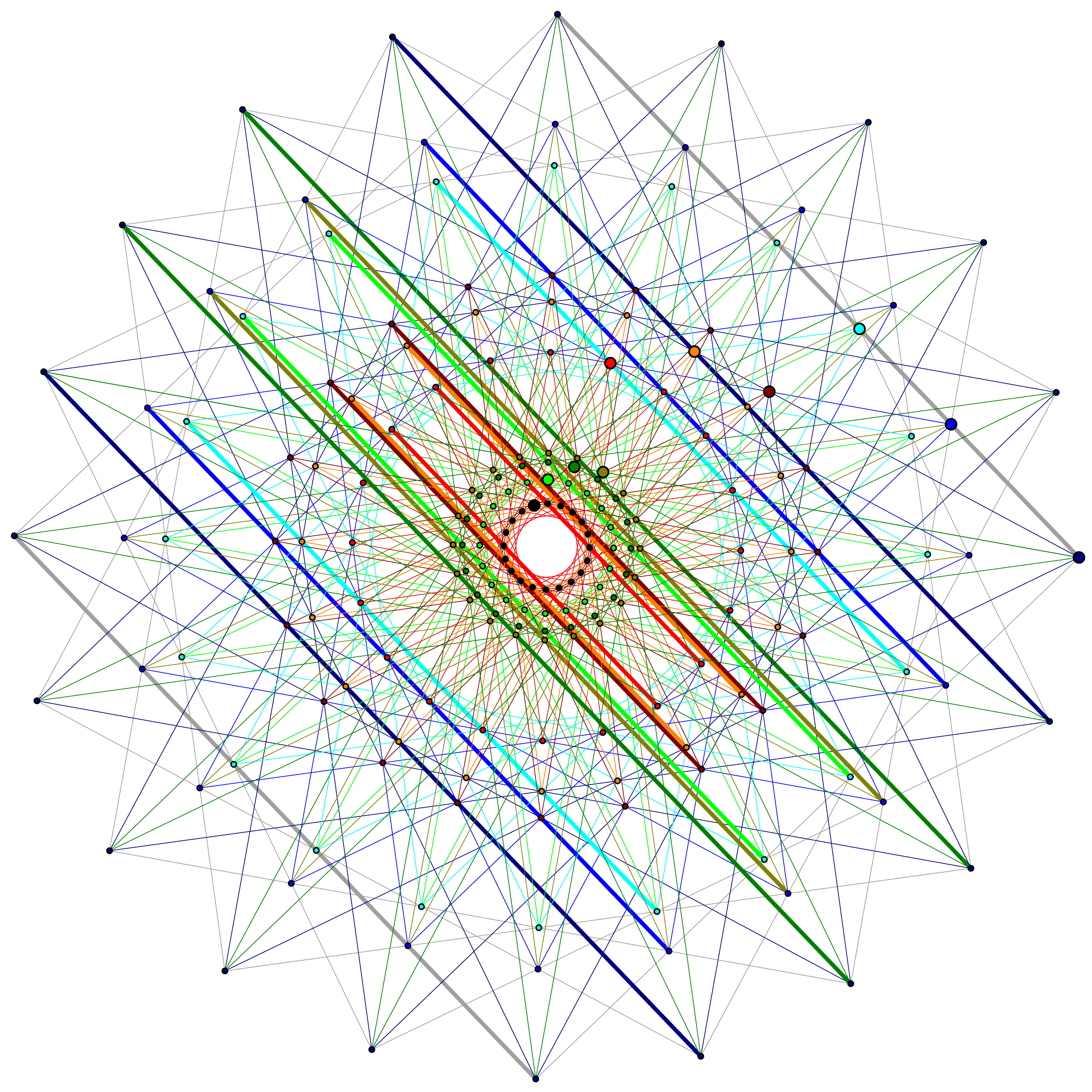}
\caption{The multicelestial configuration $20\#(5,3,1)(7,9)$ has 20 parallel lines.}
\label{default}
\end{center}
\end{figure}

\begin{thm}
$N_{6}\leq 585$.
\end{thm} 

\begin{proof}
First we enumerated the systematic 6-configurations. Then we enumerated the 5-configurations up to $1050$, following the procedure listed in the proof of Theorem \ref{thm:n5bound}. Note this count already includes the 4-configurations bootstrapped up through affine replication. We took this longer list of 5-configurations and applied affine replication to produce a starting collection of 6-configurations. 

Next, we iteratively applied the constructions \du,  $\PS$, and  $\AS$ until no new configurations were produced, to produce a ``current set'' $\mc{M}$ of missing values for $n$, a ``current set'' $\mc{S}$ of all known 6-configurations on less than 7350 points, and a subset of $\mc{S}$ called $\mc{F}$ consisting of the flexible configurations formed as outputs of affine or parallel switch or affine replication.

Next, we constructed $\du(\mc{C}, \mc{D}^{*})$ where $\mc{C} = (\ell_{6})$, $\mc{D}^{*} = (m_{6})$ such that $\ell+m-1 \in \mc{M}$, where $\mc{C} \in \mc{S}$ and $\mc{D}^{*} \in \mc{F}$.

We included these configurations into $\mc{S}$, and determined the ``current missing'' values between $n = 96$ and $n = 7350$.

Finally, we applied $\du^{(2)}(\mc C, \AR(\mc{D}))$, for $\mc{D}$ in the set of known 5-configurations, by hand, since for $n\leq 618$ we only needed to consider $(m_{5})$ 
configurations with $m\leq 100$. The resulting 6-configurations are listed in Table \ref{table:du2-6cfgs}. The configurations $\mc{C}$ listed in the table have been abbreviated by their symbol for conciseness; specifically, we used $(200_{6}) = 20\#(5,3,1)(7,9)$, $(160_{6}) = 16\#(5,4,2)(6,7)$, $(140_{6}) = 14\#(3,2,1)(4,5)$, the sporadic 4-astral $(120_{6})$, and the sporadic $(96_{6})$ configuration $\mc{A}(12;4,4;1,3,5)$, which all have halfturn symmetry.

\begin{table}[htp]
\caption{Small $6$-configurations produced by ad hoc application of $\du^{(2)}(C, \AR(D))$. }
\begin{center}
{\tiny
\begin{tabular}{l|l | l}
\hline
 $(618_{6}) = \du^{(2)}(200_{6}, \AR(60_{5}))$ &
 $(586_{6}) = \du^{(2)}(140_{6}, \AR(64_{5}))$ &
 $(566_{6}) = \du^{(2)}(120_{6}, \AR(64_{5}))$ \\
$(548_{6}) = \du^{(2)}(200_{6}, \AR(50_{5}))$&
$(542_{6}) = \du^{(2)}(96_{6},\AR(64_{5}))$&
$(536_{6}) = \du^{(2)}(160_{6}, \AR(54_{5}))$ \\
$(534_{6}) = \du^{(2)}(200_{6}, \AR(48_{5}))$ &
 $(516_{6}) = \du^{(2)}(140_{6}, \AR(54_{5}))$ &
 $(514_{6}) = \du^{(2)}(180_{6}, \AR(48_{5}))$\\
 $(494_{6}) = \du^{(2)}(160_{6}, \AR(48_{5}))$ &
 $(474_{6}) = \du^{(2)}(140_{6}, \AR(48_{5}))$ &
$(472_{6}) = \du^{(2)}(96_{6}, \AR(54_{5}))$\\ 
 $(468_{6}) = \du^{(2)}(120_{6}, \AR(50_{5}))$ &
 $(454_{6}) = \du^{(2)}(120_{6}, \AR(48_{5}))$ &
 $(444_{6}) = \du^{(2)}(96_{6}, \AR(50_{5}))$\\
\hline
\end{tabular}
}
\end{center}
\label{table:du2-6cfgs}
\end{table}%

This results in the following missing values:
$\{$97, 98, 99, 100, 101, 102, 103, 104, 105, 106, 107, 108, 109, 111, \
113, 114, 115, 116, 117, 118, 119, 121, 122, 123, 124, 125, 126, 127, \
129, 131, 132, 133, 134, 135, 136, 137, 138, 139, 141, 142, 143, 145, \
146, 147, 148, 149, 151, 152, 153, 154, 155, 156, 157, 158, 159, 161, \
162, 163, 164, 165, 166, 167, 168, 169, 171, 172, 173, 174, 175, 177, \
178, 179, 181, 182, 183, 184, 185, 186, 187, 188, 189, 193, 194, 195, \
196, 197, 198, 199, 201, 202, 203, 204, 205, 206, 207, 209, 211, 212, \
213, 214, 215, 216, 217, 218, 221, 222, 225, 226, 227, 228, 229, 231, \
232, 233, 234, 235, 236, 237, 238, 241, 242, 243, 244, 245, 246, 247, \
248, 249, 251, 252, 253, 254, 257, 258, 261, 262, 263, 264, 265, 266, \
267, 268, 269, 271, 273, 274, 275, 276, 277, 278, 281, 282, 283, 284, \
285, 289, 291, 292, 293, 294, 295, 296, 297, 298, 301, 302, 303, 305, \
306, 307, 308, 309, 311, 312, 313, 314, 315, 316, 317, 318, 321, 322, \
323, 324, 325, 326, 327, 329, 331, 332, 333, 335, 337, 338, 341, 342, \
343, 344, 345, 346, 347, 348, 349, 353, 354, 355, 356, 357, 361, 362, \
363, 364, 365, 366, 367, 369, 371, 372, 373, 374, 375, 376, 377, 385, \
386, 387, 389, 391, 393, 394, 395, 396, 397, 398, 401, 402, 403, 404, \
405, 406, 407, 408, 409, 411, 412, 413, 414, 417, 421, 422, 423, 424, \
425, 426, 427, 428, 429, 433, 434, 435, 436, 438, 441, 442, 443, \
446, 449, 451, 452, 453, 456, 457, 458, 466, 467, 471, 491, 492, 498, 502, 506, 513, 518, 522, 523, 524,  \
532, 533, 564, 584$\}. $
\end{proof}

\section{Open questions and further work}

\begin{question} What are the actual values for $N_{4}$, $N_{5}$ and $N_{6}$?\end{question}

The $\mc{A}$-series construction produces $k$-configurations for all $k \geq 4$, and the multicelestial construction can theoretically be applied to produce $2q$-configurations for $q \geq 2$, but in general, very little is known about extremely highly incident $(n_{k})$ configurations, where $k >6$.

\begin{question} What do the techniques used in this paper say about bounds on $N_{k}$ for $k \geq 7$? \end{question}

The smallest known geometric 5-configuration  has 48 points and lines, and we use it extensively in constructing additional 5- and 6-configurations. Finding a smaller configuration would likely allow lowering the value for $N_{5}$.

\begin{question} We know that the smallest geometric 3-configuration has 9 points and lines, and the smallest geometric 4-configuration has 18 points and lines, but there is no $(19_{4})$ configuration. How many points and lines does the smallest 5-configuration have?  What about the smallest $k$-configuration in general?\end{question}

\begin{question} There are only two known constructions (the $\mc{A}$-series and multicelestial constructions) that produce infinite families of symmetrically realizable 6-configurations. Can we find other symmetric constructions that will let us reduce the bound for $N_{6}$?\end{question}

\begin{question} Most of the Gr\"unbaum Calculus operations do not apply to unbalanced configurations, where the number of points and lines is different. Can we find other operations that allow the development of small bounds on the existence of $[q,k]$-configurations? \end{question}

\section{Acknowledgements} G\'abor G\'evay's research was supported in part by the Hungarian National Research,
Development and Innovation Office, OTKA grant No. SNN 132625. Toma\v{z} Pisanski's work is supported in part by the Slovenian Research Agency (research program P1-0294 and research projects N1-0032, J1-9187, J1-1690, N1-0140, J1-2481), and in part by H2020 Teaming InnoRenew CoE.

On behalf of all authors, the corresponding author states that there is no conflict of interest.

Data Availability Statement: Mathematica code and the data produced from running that code, including a list of all known $(n_{5})$ configurations with $48 \leq n \leq 576$ produced from the processes described above, is available at \url{https://github.com/UP-LaTeR/5-bounds-6-bounds}.

\bibliographystyle{plain}
\bibliography{BoundsPaper2_TO_SUBMIT}       

\section*{Appendix}

The following list gives an example of an $(n_{5})$ configuration for each $n$ where we know such a configuration exists, for $n \leq 250$.

 \begin{longtable}{l l l l}

 \caption[]
{
 Examples of known  geometric $(n_{5})$ configurations, $n \leq 250$. Values of $n$ marked with $*$ indicate that that configuration is flexible. 
For each $n$, although there may be multiple known configurations, a configuration with a highest-known number of parallel lines from our computational results is chosen to be displayed, if one is known. \label{table:n5cfgs}}\\
  \hline
   $n$ & $p$ & $q$ & Description\\ 
   \hline
 \endfirsthead
 \caption[]{Known  geometric $(n_{5})$ configurations, $n \leq 250$ (continued)}\\
   \hline
   $n$ & $p$ & $q$ & Description\\ 
 \hline
 \endhead
%
\textcolor{red}{$\leq 47$} & & &\textcolor{red}{Unknown}\\
48 & 2 & 2 & $  \mc{N}(2,6, 2)$\\ 
\textcolor{red}{$49$} & & &\textcolor{red}{Unknown}\\
50 & 1 & 1 & $  \mc{D}_4(5)$\\ 
\textcolor{red}{$51$} & & &\textcolor{red}{Unknown}\\
\textcolor{red}{$52$} & & &\textcolor{red}{Unknown}\\
\textcolor{red}{$53$} & & &\textcolor{red}{Unknown}\\
54 & 3 & 3 & $  \mc{N}'(9)$\\
\textcolor{red}{$55$} & & &\textcolor{red}{Unknown}\\ 
56 & 1 & 1 & $  \mc{A}(7)$\\ 
\textcolor{red}{$57$} & & &\textcolor{red}{Unknown}\\
\textcolor{red}{$58$} & & &\textcolor{red}{Unknown}\\
\textcolor{red}{$59$} & & &\textcolor{red}{Unknown}\\
60 & 1 & 1 & $ \mc{D}_4(6)$\\ 
\textcolor{red}{$61$} & & &\textcolor{red}{Unknown}\\
\textcolor{red}{$62$} & & &\textcolor{red}{Unknown}\\
\textcolor{red}{$63$} & & &\textcolor{red}{Unknown}\\
64 & 2 & 2 & $ \mc{A}(8)$\\ 
\textcolor{red}{$65$} & & &\textcolor{red}{Unknown}\\
66 & 1 & 1 & $ \mc{D}_5(5)$\\ 
\textcolor{red}{$67$} & & &\textcolor{red}{Unknown}\\
\textcolor{red}{$68$} & & &\textcolor{red}{Unknown}\\
\textcolor{red}{$69$} & & &\textcolor{red}{Unknown}\\
70 & 1 & 1 & $ \mc{D}_4(7)$\\ 
\textcolor{red}{$70$} & & &\textcolor{red}{Unknown}\\
72 & 2 & 2 & $ \mc{N}(2,6,3)$\\ 
\textcolor{red}{$73$} & & &\textcolor{red}{Unknown}\\
\textcolor{red}{$74$} & & &\textcolor{red}{Unknown}\\
\textcolor{red}{$75$} & & &\textcolor{red}{Unknown}\\
\textcolor{red}{$76$} & & &\textcolor{red}{Unknown}\\
\textcolor{red}{$77$} & & &\textcolor{red}{Unknown}\\
78 & 1 & 1 & $ \mc{D}_5(6)$\\ 
\textcolor{red}{$79$} & & &\textcolor{red}{Unknown}\\
80 & 2 & 2 & $ \mc{A}(10)$\\ 
81 & 3 & 3 & $ \mc{N}(3,3,3)$\\ 
\textcolor{red}{$82$} & & &\textcolor{red}{Unknown}\\
\textcolor{red}{$83$} & & &\textcolor{red}{Unknown}\\
\textcolor{red}{$84$} & & &\textcolor{red}{Unknown}\\
\textcolor{red}{$85$} & & &\textcolor{red}{Unknown}\\
\textcolor{red}{$86$} & & &\textcolor{red}{Unknown}\\
\textcolor{red}{$87$} & & &\textcolor{red}{Unknown}\\
88 & 1 & 1 & $ \mc{A}(11)$\\ 
\textcolor{red}{$89$} & & &\textcolor{red}{Unknown}\\
90 & 6 & 6 & $ \mc{N}(3,2,5)$\\ 
\textcolor{red}{$91$} & & &\textcolor{red}{Unknown}\\
\textcolor{red}{$92$} & & &\textcolor{red}{Unknown}\\
\textcolor{red}{$93$} & & &\textcolor{red}{Unknown}\\
\textcolor{red}{$94$} & & &\textcolor{red}{Unknown}\\
95 & 2 & 2 & $ \du(1)(\mc{N}(2,6,2))$\\ 
96 & 2 & 2 & $ \mc{A}(12; 5,5;1,2,4)$\\ 
\textcolor{red}{$97$} & & &\textcolor{red}{Unknown}\\
\textcolor{red}{$98$} & & &\textcolor{red}{Unknown}\\
99 & 1 & 1 & $ \du(1)(\mc{D}_4(5))$\\ 
100 & 1 & 1 & $ \mc{D}_4(10)$\\ 
\textcolor{red}{$101$} & & &\textcolor{red}{Unknown}\\
102 & 1 & 1 & $ \mc{D}_5(8)$\\ 
\textcolor{red}{$103$} & & &\textcolor{red}{Unknown}\\
104 & 1 & 1 & $ \mc{A}(13)$\\ 
\textcolor{red}{$105$} & & &\textcolor{red}{Unknown}\\
\textcolor{red}{$106$} & & &\textcolor{red}{Unknown}\\
107 & 3 & 3 & $ \du(1)(\mc{N}'(9))$\\ 
$108^*$& 18 & 0 & $ \AR(18_{4})^*$\\ 
\textcolor{red}{$109$} & & &\textcolor{red}{Unknown}\\
110 & 1 & 1 & $ \mc{D}_4(11)$\\ 
111 & 1 & 1 & $ \du(1)(\mc{A}(7))$\\ 
112 & 2 & 2 & $ \mc{A}(14)$\\ 
\textcolor{red}{$113$} & & &\textcolor{red}{Unknown}\\
114 & 1 & 1 & $ \mc{D}_5(9)$\\ 
\textcolor{red}{$115$} & & &\textcolor{red}{Unknown}\\
\textcolor{red}{$116$} & & &\textcolor{red}{Unknown}\\
\textcolor{red}{$117$} & & &\textcolor{red}{Unknown}\\
\textcolor{red}{$118$} & & &\textcolor{red}{Unknown}\\
119 & 1 & 1 & $ \du(1)(\mc{D}_4(6))$\\ 
$120^*$& 20 & 0 & $ \AR(20_{4})^*$\\ 
\textcolor{red}{$121$} & & &\textcolor{red}{Unknown}\\
\textcolor{red}{$122$} & & &\textcolor{red}{Unknown}\\
\textcolor{red}{$123$} & & &\textcolor{red}{Unknown}\\
\textcolor{red}{$124$} & & &\textcolor{red}{Unknown}\\
\textcolor{red}{$125$} & & &\textcolor{red}{Unknown}\\
$126^*$& 21 & 0 & $ \AR(21_{4})^*$\\ 
127 & 2 & 2 & $ \du(1)(\mc{A}(8))$\\ 
128 & 2 & 2 & $ \mc{A}(16)$\\ 
\textcolor{red}{$129$} & & &\textcolor{red}{Unknown}\\
130 & 1 & 1 & $ \mc{D}_4(13)$\\ 
131 & 1 & 1 & $ \du(1)(\mc{D}_5(5))$\\ 
$132^*$& 22 & 0 & $ \AR(22_{4})^*$\\ 
\textcolor{red}{$133$} & & &\textcolor{red}{Unknown}\\
\textcolor{red}{$134$} & & &\textcolor{red}{Unknown}\\
135 & 3 & 3 & $ \mc{N}(3,3,5)$\\ 
136 & 1 & 1 & $ \mc{A}(17)$\\ 
\textcolor{red}{$137$} & & &\textcolor{red}{Unknown}\\
138 & 1 & 1 & $ \mc{D}_5(11)$\\ 
139 & 1 & 1 & $ \du(1)(\mc{D}_4(7))$\\ 
140 & 1 & 1 & $ \mc{D}_4(14)$\\ 
\textcolor{red}{$141$} & & &\textcolor{red}{Unknown}\\
142 & 2 & 2 & $ \du(2)(\mc{N}(2,6,2))$\\ 
143 & 2 & 2 & $ \du(1)(\mc{N}(2,6,3))$\\ 
$144^*$& 24 & 0 & $ \AR(24_{4})^*$\\ 
\textcolor{red}{$145$} & & &\textcolor{red}{Unknown}\\
\textcolor{red}{$146$} & & &\textcolor{red}{Unknown}\\
\textcolor{red}{$147$} & & &\textcolor{red}{Unknown}\\
148 & 1 & 1 & $ \du(2)(\mc{D}_4(5))$\\ 
\textcolor{red}{$149$} & & &\textcolor{red}{Unknown}\\
$150^*$& 25 & 0 & $ \AR(25_{4})^*$\\ 
\textcolor{red}{$151$} & & &\textcolor{red}{Unknown}\\
152 & 1 & 1 & $ \mc{A}(19)$\\ 
\textcolor{red}{$153$} & & &\textcolor{red}{Unknown}\\
154 &16 & 2&$\du^{(2)}(\mc{N} (2, 6, 2), \AR(18_{4})^*)$\\
155 & 1 & 1 & $ \du(1)(\mc{D}_5(6))$\\ 
$156^*$& 26 & 0 & $ \AR(26_{4})^*$\\ 
157 & 18 & 1 & $ \du(\mc{D}_4(5),\AR(18_{4})^*)$\\ 
\textcolor{red}{$158$} & & &\textcolor{red}{Unknown}\\
159 & 2 & 2 & $ \du(1)(\mc{A}(10))$\\ 
160 & 3 & 3 & $ \du(2)(\mc{N}'(9))$\\ 
161 & 3 & 3 & $ \du(1)(\mc{N}(3,3,3))$\\ 
$162^*$& 27 & 0 & $ \AR(27_{4})^*$\\ 
163 & 18 & 1 & $ \du(\mc{A}(7),\AR(18_{4})^*)$\\ 
\textcolor{red}{$164$} & & &\textcolor{red}{Unknown}\\
\textcolor{red}{$165$} & & &\textcolor{red}{Unknown}\\
166 & 1 & 1 & $ \du(2)(\mc{A}(7))$\\ 
167 & 20 & 2 & $ \du(\mc{N}(2,6, 2),\AR(20_{4})^*)$\\ 
$168^*$& 28 & 0 & $ \AR(28_{4})^*$\\ 
169 & 20 & 1 & $ \du(\mc{D}_4(5),\AR(20_{4})^*)$\\ 
170 & 1 & 1 & $ \mc{D}_4(17)$\\ 
171 & 18 & 2 & $ \du(\mc{A}(8),\AR(18_{4})^*)$\\ 
$172$ & 19& 2&$\du^{(2)}(\mc{N}(2,6,2), \AR(21_{4})) $\\
173 & 21 & 2 & $ \du(\mc{N}(2,6, 2),\AR(21_{4})^*)$\\ 
$174^*$& 29 & 0 & $ \AR(29_{4})^*$\\ 
175 & 1 & 1 & $ \du(1)(\mc{A}(11))$\\ 
176 & 2 & 2 & $ \mc{A}(22)$\\ 
177 & 18 & 1 & $ \du(\mc{D}_4(7),\AR(18_{4})^*)$\\ 
178 & 1 & 1 & $ \du(2)(\mc{D}_4(6))$\\ 
179 & 6 & 6 & $ \du(1)(\mc{N}(3,2,5))$\\ 
$180^*$& 30 & 0 & $ \AR(30_{4})^*$\\ 
181 & 22 & 1 & $ \du(\mc{D}_4(5),\AR(22_{4})^*)$\\ 
$182$ & 18 & 2 &$\du^{(2)}(\mc{A}(8), \AR(20_{4})) $\\
183 & 20 & 2 & $ \du(\mc{A}(8),\AR(20_{4})^*)$\\ 
184 & 1 & 1 & $ \mc{A}(23)$\\ 
185 & 22 & 3 & $ \du(\mc{N}'(9), \AR(22_{4})^*)$\\ 
$186^*$& 31 & 0 & $ \AR(31_{4})^*$\\ 
187 & 22 & 1 & $ \du(\mc{A}(7),\AR(22_{4})^*)$\\ 
188 & 18 & 3 & $ \du(\mc{N}(3,3,3),\AR(18_{4})^*)$\\ 
189 & 3 & 3 & $ \mc{N}(3,3,7)$\\ 
190 & 2 & 2 & $ \du(2)(\mc{A}(8))$\\ 
191 & 2 & 2 & $ \du(1)(\mc{A}(12))$\\ 
$192^*$& 32 & 0 & $ \AR(32_{4})^*$\\ 
$193^*$& 4 & 2 & $ \AS(\mc{N}(2,6,2), 1, 0)^*$\\ 
$194^*$& 1 & 2 & $ \AS(\mc{N}(2,6,2), 2, 0)^*$\\ 
$195^*$& 1 & 4 & $ \AS(\mc{N}(2,6,2), 2, 1)^*$\\ 
$196^*$& 1 & 1 & $ \AS(\mc{N}(2,6,2), 2, 2)^*$\\ 
197 & 1 & 1 & $ \du(3)(\mc{D}_4(5))$\\ 
$198^*$& 33 & 0 & $ \AR(33_{4})^*$\\ 
199 & 1 & 1 & $ \du(1)(\mc{D}_4(10))$\\ 
200 & 1 & 1 & $ \mc{D}_4(20)$\\ 
$201^*$& 1 & 1 & $ \AS(\mc{D}_4(5), 1, 0)^*$\\ 
$202^*$& 1 & 1 & $ \AS(\mc{D}_4(5), 1, 1)^*$\\ 
203 & 1 & 1 & $ \du(1)(\mc{D}_5(8))$\\ 
$204^*$& 34 & 0 & $ \AR(34_{4})^*$\\ 
205 & 26 & 1 & $ \du(\mc{D}_4(5),\AR(26_{4})^*)$\\ 
206 & 21 & 3 & $ \du(\mc{N}(3,3,3),\AR(21_{4})^*)$\\ 
207 & 1 & 1 & $ \du(1)(\mc{A}(13))$\\ 
208 & 2 & 2 & $ \mc{A}(26)$\\ 
209 & 27 & 2 & $ \du(\mc{N}(2,6, 2),\AR(27_{4})^*)$\\ 
$210^*$& 35 & 0 & $ \AR(35_{4})^*$\\ 
211 & 27 & 1 & $ \du(\mc{D}_4(5),\AR(27_{4})^*)$\\ 
212 & 22 & 3 & $ \du(\mc{N}(3,3,3),\AR(22_{4})^*)$\\ 
213 & 3 & 3 & $ \du(3)(\mc{N}'(9))$\\ 
214 & 2 & 2 & $ \du(2)(\mc{N}(2,6,3))$\\ 
215 & 18 & 0 & $ \du(1)(\AR(18_{4})^*)$\\ 
$216^*$& 36 & 0 & $ \AR(36_{4})^*$\\ 
$217^*$& 8 & 3 & $ \AS(\mc{N}'(9), 1, 0)^*$\\ 
$218^*$& 4 & 3 & $ \AS(\mc{N}'(9), 2, 0)^*$\\ 
$219^*$& 1 & 3 & $ \AS(\mc{N}'(9), 3, 0)^*$\\ 
$220^*$& 1 & 8 & $ \AS(\mc{N}'(9), 3, 1)^*$\\ 
$221^*$& 1 & 4 & $ \AS(\mc{N}'(9), 3, 2)^*$\\ 
$222^*$& 37 & 0 & $ \AR(37_{4})^*$\\ 
223 & 2 & 2 & $ \du(1)(\mc{A}(14))$\\ 
224 & 2 & 2 & $ \mc{A}(28)$\\ 
$225^*$& 1 & 1 & $ \AS(\mc{A}(7), 1, 0)^*$\\ 
$226^*$& 1 & 1 & $ \AS(\mc{A}(7), 1, 1)^*$\\ 
227 & 1 & 1 & $ \du(1)(\mc{D}_5(9))$\\ 
$228^*$& 38 & 0 & $ \AR(38_{4})^*$\\ 
229 & 30 & 1 & $ \du(\mc{D}_4(5),\AR(30_{4})^*)$\\ 
230 & 1 & 1 & $ \mc{D}_4(23)$\\ 
231 & 28 & 2 & $ \du(\mc{A}(8),\AR(28_{4})^*)$\\ 
232 & 1 & 1 & $ \du(2)(\mc{D}_5(6))$\\ 
233 & 31 & 2 & $ \du(\mc{N}(2,6, 2),\AR(31_{4})^*)$\\ 
$234^*$& 39 & 0 & $ \AR(39_{4})^*$\\ 
235 & 31 & 1 & $ \du(\mc{D}_4(5),\AR(31_{4})^*)$\\ 
236 & 2 & 2 & $ \du(4)(\mc{N}(2,6,2))$\\ 
237 & 1 & 1 & $ \du(3)(\mc{D}_4(6))$\\ 
238 & 2 & 2 & $ \du(2)(\mc{A}(10))$\\ 
239 & 20 & 0 & $ \du(1)(\AR(20_{4})^*)$\\ 
$240^*$& 40 & 0 & $ \AR(40_{4})^*$\\ 
241 & 3 & 3 & $ \du(2)(\mc{N}(3,3,3))$\\ 
$242^*$& 1 & 1 & $ \AS(\mc{D}_4(6), 1, 1)^*$\\ 
243 & 3 & 3 & $ \mc{N}(3,3,9)$\\ 
244 & 25 & 2 & $ \du(\du(1)(\mc{N}(2,6, 2)),\AR(25_{4})^*)$\\ 
245 & 33 & 2 & $ \du(\mc{N}(2,6, 2),\AR(33_{4})^*)$\\ 
$246^*$& 41 & 0 & $ \AR(41_{4})^*$\\ 
247 & 33 & 1 & $ \du(\mc{D}_4(5),\AR(33_{4})^*)$\\ 
248 & 1 & 1 & $ \mc{A}(31)$\\ 
249 & 31 & 2 & $ \du(\mc{A}(8),\AR(31_{4})^*)$\\ 
$250^*$& 5 & 5 & $ \PS(\mc{D}_4(5))^*$\\ 
\end{longtable}

\end{document}